\documentclass[12pt]{amsart}

\usepackage{amsmath,amsfonts,amssymb}
\usepackage{color,graphicx}
\usepackage[left=2.5cm,right=2.5cm,top=3.5cm,bottom=3cm]{geometry}

\usepackage{amsthm}
\usepackage{ dsfont }
\newtheorem{theorem}{Theorem}[section]
\newtheorem{lemma}[theorem]{Lemma}
\newtheorem{definition}[theorem]{Definition}
\newtheorem{corollary}[theorem]{Corollary}
\newtheorem{proposition}[theorem]{Proposition}
\newtheorem{question}[theorem]{Question}
\newtheorem{example}[theorem]{Example}

\theoremstyle{definition}
\newtheorem{remark}[theorem]{Remark}

\usepackage[utf8]{inputenc}

\newcommand{\F}{\mathcal{F}}

\newcommand{\N}{\mathds{N}}
\newcommand{\C}{\mathds{C}}
\newcommand{\R}{\mathds{R}}
\newcommand{\Q}{\mathds{Q}}
\newcommand{\Z}{\mathds{Z}}

\newcommand{\Lip}{\textup{Lip}}

\date{\today}
\title[Recurrence and vectors escaping to infinity]{Recurrence and vectors escaping to infinity for Lipschitz operators} 
\author{Sebastián Tapia-García}

\address{ Sebastián Tapia-Garcia}
\address{Institute of Statistics and Mathematical Methods in Economics, VADOR E105-04 TU WIen, Wiedner Hauptstraße 8, A-1040 Wien.}
\email{sebastian.tapia.garcia@tuwien.ac.at}

\begin{document}

\maketitle

\begin{abstract} 
	We investigate dynamical properties of linear operators that are obtained as the linearization of Lipschitz self-maps defined on a pointed metric space. These operators are known as Lipschitz operators. 
	More concretely, for a Lipschitz operator $\widehat{f}$, we study the set of recurrent vectors and the set of vectors $\mu$ such that the sequence $(\|\widehat{f}^n(\mu)\|)_n$ goes to infinity. 
	As a consequence of our results we get that there is no wild Lipschitz operator.
	Furthermore, several examples are presented illustrating our ideas. We highlight the cases when the underlying metric space is a connected subset of $\R$ or a subset of $\mathds{Z}^d$.
	We end this paper studying some topological properties of the set of Lipschitz operators.
\end{abstract}
{\small

}
\tableofcontents
\setlength{\parindent}{0pt}

\section{Introduction}

Dynamics of linear operators is a rapidly growing area that has attracted the interest of many researchers in recent decades. 
Kitai's PhD thesis can be considered as one of the starting points of this theory, although Birkhoff and Rolewicz had studied dynamical properties of some concrete linear operators decades before, see for instance \cite{K,B,R}. 
Let $T$ be a bounded linear operator defined on a Banach space $X$ (in short $T\in \mathcal{L}(X)$) and $x\in X$. 
The orbit of $x$ under the action of $T$ is denoted by ${\textup{orb}(T,x):=\{T^nx:~n\in\N\}}$.
Cyclicity and Hypercyclicity are two of the most studied dynamical properties of linear operators:
$T\in\mathcal{L}(X)$ is called (hyper)cyclic if there is a vector $x\in X$ such that $\textup{span}(\textup{orb}(T,x))$ is dense in $X$ (resp. $\textup{orb}(T,x)$ is dense in $X$).
These two properties are intimately related to the following invariant subspace and subset problems, known as the invariant subspace problem and, respectively, the invariant subset problem: 
Does there exist an operator $T\in \mathcal{L}(X)$ such that the only closed subspaces $Y\subset X$ (resp. closed subsets) satisfying $T(Y)\subset Y$ are $\{0\}$ and $X$?
Further information about dynamics of linear operators can be found in \cite{BM,GP} and references therein.\\

In this manuscript, we are interested in the dynamical properties of the so called \textit{Lipschitz operators} (Definition~\ref{def: Lipschitz operator}). 
In order to define them, we first recall the Lipschitz-free space of a metric space $M$ (also known as Arens-Eells space or transportation cost space of $M$). 
We do not intend to provide here an exhaustive introduction to Lipschitz-free spaces. For a complete survey about Lipschitz-free spaces the reader may refer to \cite{W}.
Let $(M,d)$ and $(N,\rho)$ be two pointed metric spaces.
That is, metric spaces for which there are distinguished points denoted by $0_M\in M$ and $0_N\in N$ respectively (or just $0$ if there is no risk of confusion).
Let $L\geq0$.
A map $f:M\to N$ is called $L$-Lipschitz if it satisfies
\[\rho(f(x),f(y))\leq L d(x,y),~ \text{for all}~x,y\in M.\]
A function $f$ is called Lipschitz if there is $L\geq0$ such that $f$ is $L$-Lipschitz.
The Lipschitz constant of $f$,  denoted by $\Lip(f)$, is the least constant $L\geq 0$ for which $f$ is $L$-Lipschitz. 
The set of Lipschitz functions from $M$ to $N$ vanishing at $0$ is denoted by
\[\Lip_0(M,N):= \{f:M\to N:~ f~\text{is Lipschitz and}~f(0_M)=0_N\}.\]
If $N= \R$ (where its distinguished point is $0$), the set $\Lip_0(M,\R)$ is simply denoted by $\Lip_0(M)$. 
It is well-known that the Lipschitz constant $\Lip(\cdot)$ is a norm on $\Lip_0(M)$ for which it becomes a dual Banach space.\\

Let $(M,d)$ be a pointed metric space.
Let $\delta:M\to \Lip_0(M)^*$ be the map defined by:
\[\delta(x)(\varphi)=\varphi(x),~\forall x\in M,~\varphi\in \Lip_0(M).\]

It is well-known that $\delta$ is an (non-linear) isometry, i.e. $\|\delta(x)-\delta(y)\|=d(x,y)$ for any ${x,y\in M}$.
The Lipschitz-free space of $M$, denoted by $\F(M)$, corresponds to \[\F(M):=\overline{\textup{span}}(\delta(x):x\in M),\] where the closure is with respect to the usual norm-topology on $\Lip_0(M)^*$.
Recall that ${\delta(M\setminus\{0\})}$ is a linearly independent subset of $\F(M)$ and $\F(M)^*$ is isometrically isomorphic to $\Lip_0(M)$.
Since the Lipschitz-free space of a pointed metric space is isometrically isomorphic to the Lipschitz-free space of its completion, without loss of generality, we will assume that we always deal with complete metric spaces.
The following result is one of the most important properties of Lipschitz-free spaces:
 
\begin{theorem}\label{theo: free lipschitz}\cite{W}
	Let $(M,d)$ and $(N,\rho)$ be a pointed metric spaces and $f\in\Lip_0(M,N)$. Then, there is a unique bounded linear operator $\widehat{f}:\F(M)\to\F(N)$ such that.
	\[\delta_N(f(x))=\widehat{f}(\delta_M(x)).\]
	Moreover, $\|\widehat{f}\|=Lip(f)$.
	In particular, if $M=N$, $\widehat{f}\in\mathcal{L}(\F(M))$.
\end{theorem}

Thanks to Theorem~\ref{theo: free lipschitz}, we may think of $\F(M)$ as a linearization of the metric space $M$ and of $\widehat{f}$ as a linearization of the map $f$. 
Following the notation \cite{ACP}, we introduce the main object of our study.

\begin{definition}\label{def: Lipschitz operator}
 Let $(M,d)$ and $(N,\rho)$ be two pointed metric spaces. A bounded linear operator $T\in \mathcal{L}(\F(M),\F(N))$ is called Lipschitz operator if there is $f\in \Lip_0(M,N)$ such that $T=\widehat{f}$.
\end{definition}

Here we are interested in the relation between the dynamical properties of a given Lipschitz function $f\in \Lip_0(M,M)$ and the associated Lipschitz operator $\widehat{f}\in \mathcal{L}(\F(M))$. 
According to this line of research, properties like cyclicity, hypercyclicity, chaos, among others were recently investigated in \cite{ACP}. 
Following a similar spirit, in this paper we are mainly concerned with the set of recurrent points, the set of points with unbounded orbit and the set of points such that the sequence of iterates tends to infinity (points escaping to infinity), i.e.
\begin{align*}
	R_f&:=\{x\in M:~\liminf_{n\to \infty} d(x,f^n(x))=0\}\\
	U_f&:=\{x\in M:~\limsup_{n\to \infty} d(0,f^n(x))=\infty\}\\
	A_f&:=\{x\in M:~\lim_{n\to \infty} d(0,f^n(x))=\infty\}, 
\end{align*}
respectively.
These sets have special interest whenever $M$ is a Banach space (for instance $\F(M)$) and $f$ is a bounded linear map (for instance $\widehat{f}$). 
Indeed, let $X$ be a Banach space and $T\in \mathcal{L}(X)$.
Then, $R_T$ and $U_T$ are $G_\delta$-subsets of $X$.
Moreover, thanks to the Banach-Steinhaus theorem, $U_T$ is either empty or dense.
Also, $U_T$ and $A_T$ coincide whenever $X$ is a finite dimensional space. 
Followed by this fact, Pr\u{a}jitur\u{a} in \cite{P} proposed the following question: 
Is it true that for any $T\in\mathcal{L}(X)$, the set $A_T$ is either empty or dense in $X$? 
In 2012, Augé constructed on each infinite dimensional separable Banach space a bounded linear operator $T\in\mathcal{L}(X)$ for which the set $A_T$ is nonempty but not dense.
In fact, his operator satisfies: $X=A_T\cup R_T$ and $\textup{int}(R_T)\neq\emptyset\neq \textup{int}(A_T)$.  
Such an operator is called a \textit{wild operator}. 
Further information on wild operators can be found in \cite{A,T} and about the set $A_T$ can be found in \cite{MV,MV2}.\\

One of the first consequences of this work is the following:

\begin{theorem}\label{theo: main 1}
	There is no pointed metric space $M$ and $f\in \Lip_0(M,M)$ such that $\widehat{f}$ is a wild operator on $\F(M)$.
\end{theorem}

\begin{proof}
	See Corollary~\ref{cor: Rf dense} or Corollary~\ref{cor: examples 2}.
\end{proof}

To achieve this result, we have carried out an independent study of the set of recurrent vectors, $R_{\widehat{f}}$, and the set of vectors escaping to infinity, $A_{\widehat{f}}$. 
In concrete, as main examples we study the sets $R_{\widehat{f}}$ and $A_{\widehat{f}}$ whenever the underlying metric space is either a subset of $\mathds{Z}^d$ or $\R$.
Theorem~\ref{theo: main 2} and Theorem~\ref{theo: main 3} gather the main results of Section~\ref{sec: R set} and Section~\ref{sec: A set} respectively.

\begin{theorem}\label{theo: main 2}
	Let $(M,d)$ be a complete pointed metric space and let $f\in \Lip_0(M,M)$. Assume that $\textup{int}(R_{\widehat{f}})\neq \emptyset$. Then, $R_f= M$ and $ R_{\widehat{f}}$ is dense in $\F(M)$. Moreover, if
	\begin{enumerate}
		\item[$i)$] $M$ has at most countable many non-isolated points, or
		\item[$ii)$] $M$ is a connected subset of $\R$,
	\end{enumerate}
	then $R_{\widehat{f}}=\F(M)$. 
\end{theorem}
Further, we characterize rigid Lipschitz operators in the sense of Costakis, Manoussos and Parissis \cite{CMP}, Corollary~\ref{cor: rigid characterizetion}. 
See Section~\ref{sec: R set} for relevant definitions. 
\begin{theorem}\label{theo: main 3}
	Let $(M,d)$ be a complete pointed metric space and $f\in \Lip_0(M,M)$. Assume that $A_{\widehat{f}}\neq \emptyset$. Then, if
	\begin{enumerate}
		\item[$i)$] $M$ is a bounded space, 
		\item[$ii)$] $M$ has only bounded connected components and each ball intersects at most finitely many of them (for instance $M=\mathds{Z}^d$), or
		\item[$iii)$] $M$ is a connected subset of $\R$,

	\end{enumerate} 
	then $A_{\widehat{f}}$ is dense in $\F(M)$. Moreover,  for an arbitrary metric space $M$, if $A_{\widehat{f}}\cap\textup{span}(\delta(M))\neq\emptyset$, then $A_{\widehat{f}}$ is dense in $\F(M)$. 
\end{theorem}

Also, if $\N$ is equipped with the distance $d_\alpha$ associated to a sequence $\alpha:=(\alpha_n)_n\subset \R_+$ (defined in Section~\ref{sec: example}) and $f\in \Lip_0(\N,\N)$, then $A_{\widehat{f}}$ is either empty or dense in $\F(M)$, see Example~\ref{example 1}.\\

Finally, we study the set $\widehat{\Lip}_0(M,N):=\{\widehat{f}:~f\in\Lip_0(M,N)\}$ equipped with the weak operator topology and strong operator topology inherited from $\mathcal{L}(\F(M),\F(N))$. 
We obtain that $\widehat{\Lip}_0(M,N)$ is a WOT-closed subset of $\mathcal{L}(\F(M),\F(N))$ and that the convergences of sequences in $\widehat{\Lip}_0(M,N)$ with respect to the strong operator topology and weak operator topology coincide.\\

The structure of this paper is as follows. 
In Section~\ref{sec: example} we provide the preliminaries on Lipschitz-free spaces. Also, we present the Lipschitz-free space of $(\N,d_\alpha)$ which will be used as a source of (counter)examples.
In Section~\ref{sec: R set} we study the set of recurrent points of $\widehat{f}$, that is, $R_{\widehat{f}}$.
In Section~\ref{sec: U set} we see that in general, there is no nontrivial relation between the sets $A_f$, $U_{f}$, $A_{\widehat{f}}$ and $U_{\widehat{f}}$. 
We show that there is a metric space $(M,d)$ and $f\in \Lip_0(M,M)$ such that $U_f=M\setminus\{0\}$ but $A_{\widehat{f}}=\emptyset$ (even more, $U_{\widehat{f}}=\F(M)\setminus \{0\}$).
Also, by considering a bounded metric space, we show that there is a metric space $(M,d)$ and $f\in\Lip_0(M,M)$ such that $U_f=\emptyset$ but $A_{\widehat{f}}$ is dense in $\F(M)$.
In Section~\ref{sec: A set} we study the set $A_{\widehat{f}}$.
In Section~\ref{sec: top} we study the set $\widehat{\Lip}_0(M,N)$ as a subset of $\mathcal{L}(\F(M),\F(N))$, considering the strong operator topology and the weak operator topology.
We end this manuscript by stating some perspectives which arise from this work.\\

\textbf{Notation:} $(M,d)$ always denotes a complete pointed metric space and $\F(M)$ denotes its Lipschitz-free space.
The symbols $x,y,z$ will denote points in $M$; $u,v$ vectors in $\textup{span}(\delta(M))$; $\mu,\nu$ arbitrary vectors in $\F(M)$.
The symbols $f,g$ will always denote maps in $\Lip_0(M,M)$ or $\Lip_0(M,N)$, where $N$ is another complete pointed metric space. 
The symbol $\varphi$ will denote an element in $\Lip_0(M)$.
The symbol $~\widehat{\empty}~{\empty}~$denotes the functor from $\Lip_0(M,N)$ to $\mathcal{L}(\F(M),\F(N))$. The symbol $\delta$ denotes the canonical isometry from $M$ to $\F(M)$.
Let $x\in M$ and $A\subset M$ be a nonempty set. 
By $\overline{A}$ and $\textup{int}(A)$ we denote the closure and interior of $A$ respectively.
We write $\textup{dist}(x,A):=\inf_{y\in A}d(x,y)$ and $\textup{diam}(A)=\sup\{d(y,z):~y,z\in A\}$.
For $r>0$, the open and closed ball centered at $x$ of radius $r$ are denoted by
$B(x,r)$ and $\overline{B}(x,r)$ respectively, i.e. $\{y\in M:~ d(x,y)<r\}$ and $\{y\in M:~ d(x,y)\leq r\}$ respectively. 
We denote by $\N$ the set of positive integers starting from $0$. 
$\R_+$ denotes the set of positive real numbers. 
For two Banach spaces $X$ and $Y$, we denote by $\mathcal{L}(X,Y)$ the space of bounded linear operators defined from $X$ to $Y$ and by $\mathcal{L}(X)$ the space $\mathcal{L}(X,X)$.\\\

\section{Preliminary on Lipschitz-free spaces}\label{sec: example}

Let $(M,d)$ be a complete pointed metric space with distinguished point $0$.
As usual, for a Lipschitz function $\varphi\in \Lip_0(M)$, its support is defined by
$\textup{supp}(\varphi)=\overline{\{x\in M:~\varphi(x)\neq 0\}}$.
A recent work developed a useful notion of support for a vector $\mu\in \F(M)$.

\begin{definition}\cite[Definition 2.5]{APPP}
	Let $\mu\in \F(M)$. The support of $\mu$ is defined by
	\[\textup{supp}(\mu):= \bigcap\{K:~ K\subset M~\text{closed},~\mu\in \overline{\textup{span}}(\delta(x):~x\in K)\}.\]
\end{definition}

One of the main properties of this definition of support is the following \cite[below Definition 2.5]{APPP}: if $\mu\in \F(M)\setminus\{0\}$, then \[\mu\in\overline{\textup{span}}(\delta(x):~x\in \textup{supp}(\mu)).\]
Let us provide some examples. 
If $\mu=0$, then $\textup{supp}(\mu)=\emptyset$. 
If $\mu\in \textup{span}(\delta(M))\setminus\{0\}$, there are a one-to-one finite sequence $(x_i)_{i=1}^k\subset M\setminus\{0\}$ and $(\lambda_i)_{i=1}^k\in \R\setminus\{0\}$ such that $\mu=\sum_{i=1}^k \lambda_i\delta(x_i)$.
In this case, $\textup{supp}(\mu)=\{x_i:~i=1,...,k\}$.
We need the following proposition.
\begin{proposition}\cite[Proposition 2.7]{APPP}\label{prop: supp mu}
Let $\mu\in \F(M)$. 
Then, $x\in \textup{supp}(\mu)$ if and only if, for any $\varepsilon>0$ there is $\varphi\in\Lip_0(M)$, such that $\textup{supp}(\varphi)\subset B(x,\varepsilon)$ and $\langle \varphi,\mu\rangle \neq 0$.  
\end{proposition}

Now, we present a family of metric spaces for which its Lipschitz-free space is well understood. 
In the forthcoming sections we will use different kind of functions defined on these metric spaces, with either recurrent properties or (forward, backward) shift behavior. \\

Let $\alpha=(\alpha_n)_n\subset\R_+$. We define the distance $d_\alpha$ on $\mathds{N}$ by 
\[d_\alpha(i,j)=
\begin{cases}
	0&~\text{if } i=j\\
	\alpha_i&~\text{if } j=0,~i\neq 0\\
	\alpha_j&~\text{if } i=0,~j\neq 0\\
	\alpha_i+\alpha_j&~\text{if }  0\neq i\neq j\neq 0.
\end{cases}
\]

The following propositions gather some simple facts about the space $(\mathds{N},d_\alpha)$. 

\begin{proposition}\label{prop: function on N alpha}
	Let $(\alpha_n)_n\subset \R_+$, $f:(\N,d_\alpha)\to(\N,d_\alpha)$ and $\varphi:(\N,d_\alpha)\to \R$ such that $f(0)=0$ and $\varphi(0)=0$. Then,
	\begin{enumerate}
		\item[$i)$] $f$ is $L$-Lipschitz if and only if $\alpha_{f(n)}\leq L\alpha_n$ for all $n\in \N$ 
		\item[$ii)$] $\varphi$ is $L$-Lipschitz if and only if $|\varphi(n)|\leq L\alpha_n$ for all $n\in \N$.
	\end{enumerate}
\end{proposition}
\begin{proof}
	$i$): If $f$ is $L$-Lipschitz, then $\alpha_{f(n)}=d(f(n),0)\leq Ld(n,0)=L\alpha_n$. 
	Conversely, if $0\neq n\neq m\neq 0$ and assuming $\alpha_0=0$ if necessary, we have that \[d(f(n),f(m))\leq \alpha_{f(n)}+\alpha_{f(m)}\leq L(\alpha_n+\alpha_m)=Ld(n,m).\]
	Now, it readily follows that $f$ is an $L$-Lipschitz function.\\
	$ii$): It is analogous to (1).
\end{proof}

The following proposition is a straightforward consequence of \cite[Proposition 3.9]{W} or more precisely \cite[Proposition 1.6]{ACP}.

\begin{proposition}\label{free of N}
	Let $\alpha=(\alpha_n)_n\subset \R_+$. 
	Then
	\begin{enumerate}
		\item[$i)$] $\mathcal{F}(\mathds{N},d_\alpha)$ is isometrically isomorphic to $\ell^1$, and
		\item[$ii)$] for any $\mu\in \mathcal{F}(\mathds{N},d_\alpha)$, there is $\lambda\in \ell^1$ such that 
		\[\mu= \sum_{n=1}^\infty \lambda_n \frac{\delta(n)}{\alpha_n}~\text{and }\|\mu\| = \sum_{n=1}^{\infty} |\lambda_n|.\]
		In other words, $\delta(n)$ can be identified with $\alpha_ne_n\in\ell^1$ for all $n\geq 1$, where $(e_n)_n$ denotes the canonical basis of $\ell^1$.
	\end{enumerate}
\end{proposition}

\section{Recurrent vectors}\label{sec: R set}

In this section we are interested in the set of recurrent vectors for Lipschitz self maps $f$. We focus our study in those functions $f$ such that $\textup{int}(R_{\widehat{f}})\neq \emptyset$. 
Let us start with the definitions of recurrent and rigid operator. 
Further information about these kind of operators can be found in \cite{CMP}.

\begin{definition}
	Let $X$ be a Banach space and $T\in \mathcal{L}(X)$. We say that $T$ is recurrent if $R_T$ is dense in $X$. We say that $T$ is rigid if there is an increasing sequence $(n(j))_j\subset\N$ such that $(T^{n(j)}x)_j$ converges to $x$ for all $x\in X$, i.e. $(T^{n(j)})_j$ converges to the identity operator in the strong operator topology.
\end{definition}

Our first goal is to show that if $\textup{int}(R_{\widehat{f}})\neq \emptyset$, then $R_f=M$ and $\widehat{f}$ is a recurrent operator. In order to do this, we need the following proposition.

\begin{proposition}\label{prop: recurrent points interior}
	Let $(M,d)$ be a pointed metric space and $f\in \Lip_0(M,M)$ be such that $R_{\widehat{f}}$ has nonempty interior. 
	Then, for every nonempty and finite set $N\subset M$, there is an increasing sequence $(n(j))_j\subset\N$ such that $(f^{n(j)}(x))_j$ converges to $x$ for all $x\in N$. 
	Particularly, $R_f=M$.
\end{proposition}

\begin{proof}
Since $R_{\widehat{f}}$ has nonempty interior and $\textup{span}(\delta(M))$ is dense in $\F(M)$, there is $u=\sum_{i=1}^k\lambda_i\delta(x_i)\in \text{int}(R_{\widehat{f}})$ where $(x_i)_{i=1}^k$ is a one-to-one finite sequence in $M\setminus\{0\}$ and $(\lambda_i)_{i=1}^k\subset\R\setminus\{0\}$ is $\mathds{Q}$-linearly independent.
Thus, we have that $\sum_{i\in I}\lambda_i\notin \{\lambda_j:~j\notin I\}$ for any $I\subset\{1,...,k\}$, where we use the convention $\sum_{i\in\emptyset}\lambda_i=0$.\\

\textup{Claim:} There is an increasing sequence $(n(j))_j\subset\N$ such that $(f^{n(j)}(x_i))_j$ converges to $x_i$ for all $i=1,...,k$.\\

Indeed, since $u\in R_{\widehat{f}}$, there is an increasing sequence $(n(j))_j\subset \N$ such that $(\widehat{f}^{n(j)}(u))_j$ converges to $u$. 
Let us see that $(f^{n(j)}(x_i))_j$ converges to $x_i$ for all $i=1,...,k$.
Assume, towards a contradiction, that there is $i\in\{1,...,k\}$ such that $(f^{n(j)}(x_i))_j$ does not converge to $x_i$.
Without loss of generality, consider $i=1$. 
Thus, there are a subsequence of $(n(j))_j$, still denoted by $(n(j))_j$, and $\varepsilon>0$ such that $d(f^{n(j)}(x_1),x_1)>\varepsilon$ for all $j\in\N$. 
A straightforward induction gives us that, there are $\varepsilon>0$ and a subsequence of $(n(j))_j$, still denoted by $(n(j))_j$, such that, for any $i\geq 2$,

\[\lim_{j\to\infty}f^{n(j)}(x_i)= x_1,~\text{or } d(f^{n(j)}(x_i),x_1)>\varepsilon,~\forall ~j\in\N,\]
and $d(f^{n(j)}(x_1),x_1)>\varepsilon$ for all $j\in\N$.
Let $\widetilde{\varepsilon}=\min(d(x_i,x_j):~1\leq i\neq j\leq k)$, and $\varphi\in \textup{Lip}_0(M)$ such that $\varphi(x_1)=1$ and $\textup{supp}(\varphi)\subset B(x_1,\min(\varepsilon,\widetilde{\varepsilon})/2)$. Let $I\subset \{1,...,k\}$ be the set indexes $i$ such that $(f^{n(j)}(x_i))_j$ converges to $x_1$. Thus
\[\lambda_1=\lim_{j\to\infty} \langle \varphi, \widehat{f}^{n(j)}(u)\rangle =\lim_{j\to\infty} \sum_{i=1}^k \lambda_i \varphi(f^{n(j)}(x_i))= \sum_{i\in I} \lambda_i \neq \lambda_1. \]
In the last part of above expression we use that $1\notin I$. This ends the proof of the claim.\\

Let us finish the proof of the proposition. 
Let $N\subset M\setminus\{0\}$ be a finite subset. 
Let $N'=N\setminus\{x_i:i=1,...,k\}$ and $\{\alpha_x:~x\in N'\}\subset \R$ such that $\{ \lambda_i:~i=1,...,k\}\cup\{\alpha_x:~x\in N'\}$ is $\mathds{Q}$-linearly independent. 
Since $u\in \text{int}(R_{\widehat{f}})$, there is $p\in \N$ such that 
\[u+\dfrac{1}{p}\sum_{x\in N'}\alpha_x\delta(x)\in R_{\widehat{f}}.\]
Now, the claim gives us a sequence $(n(j))_j\subset \N$ such that $(f^{n(j)}(x))_j$ converges to $x$ for all $x\in N$.
\end{proof}

As a consequence we get the following result. 
\begin{corollary}~\label{cor: Rf dense}
 Let $(M,d)$ be a pointed metric space and $f\in\Lip_0(M,M)$. If $R_{\widehat{f}}$ has nonempty interior, then $\textup{span}(\delta(M))\subset R_{\widehat{f}}$. Thus, $\widehat{f}$ is a recurrent operator. Particularly, there is no wild operator of the form $\widehat{f}$.
\end{corollary}

\begin{proof}
If $\textup{int}(R_{\widehat{f}})$ is nonempty, Proposition~\ref{prop: recurrent points interior} readily implies that $\textup{span}(\delta(M))\subset R_{\widehat{f}}$. Hence, $R_{\widehat{f}}$ is a dense subset of $\F(M)$.
Moreover, since $A_{\widehat{f}}\subset \F(M)\setminus R_{\widehat{f}}$, $A_{\widehat{f}}$ has empty interior. So $\widehat{f}$ is not a wild operator.
\end{proof}

Before continuing, let us recall the strong and weak operator topologies on the space of linear operators.
\begin{definition}\label{def: sot and wot}
	Let $X$ and $Y$ be two Banach spaces, $(T_n)_n\subset \mathcal{L}(X,Y)$ and $T\in \mathcal{L}(X,Y)$. 
	We say that $(T_n)_n$ converges to $T$ with respect to the strong operator topology (SOT) if for any $x\in X$, $(T_nx)_n$ converges to $Tx$.
	We say that $(T_n)_n$ converges to $T$ with respect to the weak operator topology (WOT) if for any $x\in X$ and $y^*\in Y^*$, $(\langle y^*, T_nx\rangle)_n$ converges to $\langle y^*,Tx\rangle$.
\end{definition}
Note that, in general, any SOT-convergent sequence is also WOT-convergent. In what follows we characterize rigid operators.
\begin{proposition}\label{prop: bounded Lip constant}
	Let $(M,d)$ and $(N,\rho)$ be complete pointed metric spaces, $(f_n)_n\subset \Lip_0(M,N)$, $D\subset M$ be a dense subset and $g:D\to N$ be a function. 
	If $(f_{n}(x))_n$ converges to $g(x) $ for all $x\in D$ and there is $C>0$ such that $\Lip(f_{n})\leq C$ for all $j$, then $(\widehat{f}^{n})_n$ is SOT-convergent.
\end{proposition}

\begin{proof}
	Since $g:D\to N$ is the pointwise limit of a bounded sequence of Lipschitz functions, $g$ is in fact a Lipschitz function, with $\Lip(g)\leq C$. 
	Moreover, since $N$ is a complete metric space, $g$ can be uniquely extended to $\overline{D}=M$, and the extension keep the same Lipschitz constant. 
	Let us denote this extension by $g$.
	We claim that $f_{n}(y)\to g(y)$ for all $y\in M$. 
	Indeed, fix $y\in M$. 
	Note that 
	\begin{align*}
	\rho(f_{n}(y),g(y))&\leq \rho(f_{n}(y),f_{n}(x))+\rho(f_{n}(x),g(x))+\rho(g(x),g(y))\\
	&\leq \rho(f_{n}(x),g(x))+2Cd(x,y)	
	\end{align*}

	Thus, $\limsup_n \rho(f_{n}(y),y)\leq 2Cd(x,y)$, for all $x\in D$.
	Since $D$ is dense in $M$, the claim follows.\\
	
	Now we show that $(\widehat{f}_n)_n$ converges to $\widehat{g}$ for the strong operator topology.
	Let $\mu\in \F(M)$. Since $\textup{span}(\delta(M))$ is a dense subspace of $\F(M)$, there is sequence $(u_k)_k\subset \textup{span}(\delta(M))$ satisfying 
	\[ \| u_k\| \leq \dfrac{\|\mu\|}{2^{k}}~\text{for all }k\in\N,~\text{and } \mu=\sum_{k=0}^\infty u_k. \]
	Note that the series is absolutely convergent. Thus, $T(\mu)=\sum_{k=0}^\infty T(u_k)$ for any $T\in\mathcal{L}(\F(M))$.
	Also, for any $k\in\N$, since $u_k\in \textup{span}(\delta(M))$, it follows that $(\widehat{f}_{n}(u_k))_n$ converges to $\widehat{g}(u_k)$. 
	Let $l\geq 1$.
	We have that
	\begin{align*}
	\limsup_{n\to\infty}\| \widehat{f}_{n}(\mu)-\widehat{g}(\mu)\|&\leq \limsup_{n\to\infty}\|\sum_{k=0}^{l} \widehat{f}_{n}(u_k)-\widehat{g}(u_k)\| + 2C\sum_{k=l+1}^\infty \|u_k\|\leq \dfrac{2C}{2^{l}}\|\mu\|.
	\end{align*}
	Since $l$ is arbitrary, we get that $(\widehat{f}_{n}(\mu))_j$ norm-converges to $\widehat{g}(\mu)$.
	Thus, $(\widehat{f}_n)_n$ SOT-converges to $\widehat{g}$.
\end{proof}

\begin{remark}\label{rem: prop bounded lip}
	The following converse of Proposition~\ref{prop: bounded Lip constant} also holds: if $g\in\Lip_0(M,N)$ and $(\widehat{f}_{n})_n$ SOT-converges to $\widehat{g}\in \mathcal{L}(\F(M))$, then there is $C>0$ such that $\|\widehat{f}_n\|<C$ (and thus $\Lip(f_n)<C$) for all $n$ and $(f_n(x))_n$ converges to $g(x)$ for all $x\in M$. 
\end{remark}
Now, let us apply Proposition~\ref{prop: bounded Lip constant} to get rigid operators.

\begin{corollary}\label{cor: rigid characterizetion}
	Let $(M,d)$ be a complete pointed metric space and $f\in \Lip_0(M,M)$. 
	Then, $\widehat{f}$ is a rigid operator if and only if there are an increasing sequence $(n(k))_k\subset \N$ and a constant $C\geq 0$ such that $\Lip(f^{n(k)})\leq C$ for all $k\in \N$ and $(f^{n(k)}(x))_k$ converges to $x$ for all $x\in M$.    
\end{corollary}

\begin{proof}
	If $\widehat{f}$ is a rigid operator, then there is a sequence $(n(k))\subset \N$ such that $(\widehat{f}^{n(k)}(x))_k$ SOT-converges to the identity $Id$. This readily implies the existence of the constant $C\geq 0$ such that $\Lip(f^{n(k)})=\|\widehat{f}^{n(k)}\|\leq C$ for all $k\in \N$ and that $(f^{n(k)}(x))_k$ converges to $x$ for all $x\in M$.
	On the other hand, let us assume that there are $C\geq 0$ and $(n(k))_k\subset\N$ as in the statement of the proposition. 
	We show that $\widehat{f}$ is a rigid operator. 
	Indeed, it readily follows by using Proposition~\ref{prop: bounded Lip constant}. 
	It is enough to use the sequence $(f^{n(k)})_k$, $D=M$ and $g=id$.
\end{proof}

\begin{corollary}\label{cor: rigid} 
	Let $(M,d)$ be a pointed metric space and $f\in \Lip_0(M,M)$ such that $\textup{int}(R_{\widehat{f}})\neq \emptyset$. Assume that there is $C>0$ such that $\Lip(f^n)\leq C<\infty$ for all $n\in \N$. Then, $R_{\widehat{f}}=\F(M)$.
	Moreover, if $M$ is separable, $\widehat{f}$ is a rigid operator.
\end{corollary}

\begin{proof}
	Thanks to Proposition~\ref{prop: recurrent points interior}, we know that for any finite subset $N\subset M$, there is a sequence $(n(j))_j$ such that $(f^{n(j)}(x))_j$ converges to $x$ for all $x\in N$. 
	A standard diagonal-like argument shows that for any countable subset $D\subset M$ there is a sequence $(n(j))_j$ such that $(f^{n(j)}(x))_j$ converges to $x$ for all $x\in D$.
	Let $\mu\in \F(M)$ be an arbitrary vector in the Lipschitz free space of $M$.
	As in the proof of Proposition~\ref{prop: bounded Lip constant}, there is $(u_k)_k\subset \textup{span}(\delta(M))$  such that $\|u_k\| \leq 2^{-k}\|\mu\|$ for all $k\in\N$ and
	$\mu=\sum_{k=0}^\infty u_k$. 
	For each $k\geq 0$, let $I_k=\textup{supp}(u_k)$ which is a finite set. 
	Let us now consider $D=\cup_{k=0}^{\infty}I_k$ and $(n(j))_j\subset \N$ be a sequence such that $(f^{n(j)}(x))_j$ converges to $x$ for all $x\in D$.
	Following the same argument of the proof of Proposition~\ref{prop: bounded Lip constant} (replacing $g(\mu)$ by $\mu$), we get that $(\widehat{f}^{n(j)}(\mu))_j$ converges to $\mu$.
	Therefore, $\mu\in R_{\widehat{f}}$.
	Finally, if $M$ is separable, $D$ can be chosen as a countable dense subset of $M$ and we fix the sequence $(n(k))_k$ as above. Then, the conclusion follows directly from Proposition~\ref{prop: bounded Lip constant} by considering the sequence $(f^{n(k)})_k$ and $g=id$.
\end{proof}
\begin{remark} Let us check that the assumptions of Proposition~\ref{prop: bounded Lip constant} do not imply the convergence in norm of the sequence $(\widehat{f}_n)_n$. 
Indeed, let us consider $M=[0,1]$ and $f_n=\textup{dist}(\cdot,[0,n^{-1}])\in \Lip_0(M,M)$. It follows that $(f_n)_n$ converges to the identity pointwise. Also, $\textup{Lip}(f_n)=1$ for all $n\geq 1$. 
However, $(\widehat{f}_n)_n$ does not converge to $Id\in\mathcal{L}(\F(M))$ with respect the operator norm. 
Indeed,
\[\|\widehat{f}_n-Id\| \geq \frac{\|(\widehat{f}_n-Id)(\delta((2n)^{-1}))\|}{(2n)^{-1}}=1,~\text{for all }n\geq 1 .\]
\end{remark} 
The following example shows that the hypothesis of uniform boundedness of the Lipschitz constant can not be removed from Corollary~\ref{cor: rigid}.

\begin{example}\label{example 3}
	There is a pointed metric space $(M,d)$ and a function $f\in \Lip_0(M,M)$ such that there is an increasing sequence $(n(j))_j\subset \N$ satisfying $f^{n(j)}(x)\to x$ for all $x\in M$ but $R_{\widehat{f}}\neq \mathcal{F}(M)$.
\end{example}

\begin{proof}
	Let us use the metric space defined in Section~\ref{sec: example}, $(\N,d_\alpha)$. Let $(s_n)_n$ be the sequence defined by $s_n=n(n+1)/2$, that is, $s_1=1$ and $s_{n+1}-s_n=n+1$ for all $n\geq 1$.
	Let us define the sequence $\alpha=(\alpha_n)_n\subset \R_+$.
	Let $m\in \N$, with $m\geq 1$, and $n\in \N$ such that $m\in [s_n,s_{n+1})$. Then, $\alpha_m=2^{m-s_n}$. Let $f:\N\to\N$ defined by

	\[f(m):=\begin{cases}
	0&~\text{if } m=0,\\
	s_n & ~ \text{if } \exists n\in\N,~m=s_{n+1}-1,\\
	m+1& ~\text{otherwise.}
	\end{cases}\]

	Observe that $ f^{s_n}(m)=m$ for all $m\in [s_n,s_{n+1})$. Thus, every point in $\N$ is periodic for $f$.  
	Also, thanks to Proposition~\ref{prop: function on N alpha}, $\Lip(f)=2$. 
	Moreover, it also follows that $\Lip(f^n)=2^n$ for every $n\in \N$, and thus, $\|\widehat{f}^n\|=2^n$. 
	Now, due to \cite[Theorem 3]{MV}, $A_{\widehat{f}}$ is dense in $\F(M)$.
\end{proof}

\begin{remark}Despite the fact that we use a general result of Müller and Vr\u{s}ovsk\'{y} to get the desired conclusion, we can explicit a vector in $A_{\widehat{f}}$, namely $\mu=\sum_{n=1}^\infty \delta(s_n)/n^2$. Moreover, since every orbit generated by $f$ is bounded (in fact finite), it follows that $\mu+\textup{span}(\delta(\N))\subset A_{\widehat{f}}$. Thus, $A_{\widehat{f}}$ is dense in $\F(M)$. For details see Proposition~\ref{bounded orbits}.
\end{remark}

The following example shows that the hypothesis of separability can not be removed from Corollary~\ref{cor: rigid}.

\begin{example}\label{example nonseparable}
	There is a nonseparable pointed metric space $M$, $f\in \Lip_0(M,M)$ such that $R_{\widehat{f}}=\F(M)$ and $\Lip(f^n)=1$ for all $n\in \N$, but $\widehat{f}$ is not a rigid operator.
\end{example}

We need the following two lemma in order to construct our example.
\begin{lemma}\label{lem: common sequence uncountable}
	Let $(\alpha_p)_{p=1}^q\subset [0,2\pi)$. Then, there is an increasing sequence $(n(k))_k\subset \N$ such that $(\exp(n(k)\alpha_p i))_k$ converges to $1$ for all $p=1,...,q$.
\end{lemma}
\begin{proof}If $\alpha_p\in \Q\pi$ for all $p\in\{1,...,q\}$, then the lemma easily follows. 
Let $A$ be a maximal subset of $\{\alpha_p:p=1,...,q\}$ such that $A\cup\{\pi\}$ is $\Q$-linearly independent. 
Assume that $A$ is a nonempty set.
Up to rearrangement of $(\alpha_p)_{p=1}^q$, let $r=|A|$ and assume that $A=\{\alpha_p:p=1,...r\}$.
Let $(a_{p,j})_{p=r+1,~j=1}^{q, ~\hspace{0.7cm}r+1}\subset \Q$ be such that
\[\alpha_p = \sum_{j=1}^r a_{p,j}\alpha_j + a_{p,r+1}\pi,~\text{for all } p=r+1,...,q.\]
Let $N\in \N$ be such that $Na_{p,j}\in 2\Z$ for all $p,j$.
Now, thanks to Kronecker's Theorem, there is a sequence $(m(k))_k$ such that 
$(\exp(m(k)\alpha_p i ))_k$ converges to $1$ for all $p=1,...,r$. 
Set $n(k)=Nm(k)$ for all $k\in \N$. It follows that $(\exp(n(k)\alpha_p i ))_k$ converges to $1$ for all $p=1,...,r$ and also that
\[\exp( n(k)\alpha_p i) = \prod_{j=1}^r \exp (Na_{p,j} \alpha_j m(k) i)= \prod_{j=1}^r \exp (\alpha_j m(k) i)^{(Na_{p,j})}, \]
tends to $1$ as $k$ tends to $\infty$ for all $p=r+1,...,q$. 
\end{proof}
\begin{lemma}\label{lem: nonconvergent uncountable}
	For any increasing sequence $(n(k))_k\subset \N$, there is $\alpha\in [0,2\pi)$ such that $(\exp(n(k)\alpha i))_k$ does not converge to $1$.
\end{lemma}
\begin{proof}
	Let us proceed towards a contradiction.
	Assume that there is an increasing sequence $(n(k))_k\subset \N$ such that $(\exp(n(k)\alpha i))_k$ converges to $1$ for all $\alpha\in [0,2\pi)$.
	Let $N\in \N$ be a positive integer. 
	Note that if we consider $\alpha= 2\pi/N$, we get that there is $K\in \N$ such that $n(k)\in N\N$ for all $k\geq K$. 
	Let $(m(k))_k$ be a subsequence of $(n(k))_k$ such that for all $k\geq 1$:
	\begin{align}\label{properties m(k)}
	2m(k)| m(k+1),~\hspace{0.5cm} \text{and } \hspace{0.5cm} m(k)\sum_{j=k+1}^\infty \dfrac{1}{m(j)}< 2^{-k}.	
	\end{align}
	Consider now 
	\[\widetilde{\alpha}:=2\pi \sum_{j=1}^\infty \dfrac{3+(-1)^j}{4m(j)}. \]
	Note that, for any $k\geq 1$ 
	\begin{align*}
		m(k)\widetilde{\alpha}&= 2\pi \left( m(k)\sum_{j=1}^{k-1} \dfrac{3+(-1)^j}{4m(j)} + \dfrac{3+(-1)^k}{4} + m(k)\sum_{j=k+1}^{\infty}\dfrac{3+(-1)^j}{4m(j)} \right)\\
		&= 2\pi ( z(k) + 4^{-1}(3+(-1)^k) + \varepsilon(k)),
	\end{align*}
	where $z(k)$ and $\varepsilon(k)$ correspond to the first and third term of the right hand side of the above line respectively.
	Thanks to \eqref{properties m(k)}, we have that $z(k)\in \N$ for all $k\geq 1$ and $(\varepsilon(k))_k$ converges to $0$. 
	Therefore, the sequence $(\exp(m(k)\widetilde{\alpha}i))_k$ has two accumulation points, namely $1$ and $-1$. 
	Thus $(\exp(n(k)\widetilde{\alpha}i))_k$ is not convergent, which is a contradiction.
\end{proof}
Now we can present our Example~\ref{example nonseparable}.
\begin{proof}[Proof of Example~\ref{example nonseparable}]
	Let $S$ be the unit circumference of $\C$, i.e. $S=\{e^{i\theta}:~\theta\in[0,2\pi)\}$.
	Let $M$ be the disjoint union of $S\times[0,2\pi)$ and $\{0\}$, i.e. $M:= S\times[0,2\pi) \dot{\cup} \{0\}$. Let us define the metric $d$ on $M$ by: 
	\[d(x,y)=\begin{cases}
		0&~\text{if }x=y\\
	 	1&~\text{if }x=0\neq y\\
	 	1&~\text{if }y=0\neq x\\
	 	|e^{i\theta_1}-e^{i\theta_2}|&~\text{if } x=(e^{i\theta_1},\alpha),~ y=(e^{i\theta_2},\alpha)\\
	 	2 &~\text{if } x=(e^{i\theta_1},\alpha),~y=(e^{i\theta_2},\beta),~ \alpha\neq \beta.
	\end{cases}\] 
It easily follows that $(M,d)$ is a complete, nonseparable pointed metric space ($0_M=0$). 
In fact, $M$ can be seen as the disjoint union of uncountable many unit circumferences.
Let $f\in \Lip_0(M,M)$ be the function defined by $f(0)=0$ and $f((e^{\theta i},\alpha))=(\exp(i(\theta+\alpha)),\alpha)$ for any $\theta\in \R$ and $\alpha\in[0,2\pi)$. 
Note that $f(S\times\{\alpha\})=S\times\{\alpha\}$. 
It follows that $\Lip(f^n)=1$ for all $n\geq 1$ and $R_f=M$. 
Using Lemma~\ref{lem: common sequence uncountable} and a straightforward modification of the second part of the proof of Proposition~\ref{prop: bounded Lip constant} we get that $R_{\widehat{f}}=\F(M)$. However, thanks to Lemma~\ref{lem: nonconvergent uncountable}, $\widehat{f}$ is not rigid because there is no an increasing sequence $(n(k))_k\subset \N$ such that $(\widehat{f}^{n(k)}(\delta(x)))_k$ converges to $\delta(x)$ for all $x\in M$. 
\end{proof}
Thanks to Corollary~\ref{cor: rigid characterizetion}, if $\widehat{f}$ is rigid then $R_f=M$. The following example shows that this is not a property that the recurrent Lipschitz operators satisfy.
\begin{example}
	There is a pointed metric space $M$ and $f\in \Lip_0(M,M)$ such that $R_f=\{0\}\neq M$ and $\widehat{f}$ is recurrent. Moreover, $\widehat{f}$ is hypercyclic.
\end{example}

\begin{proof}
	Consider $M=\{0\}\cup\{2^{-n}:~n\in \N\}$ equipped with the metric inherited from $\R$. 
	Let $f:M\to M$ defined by $f(0)=f(1)=0$ and $f(2^{-n})=f(2^{-n+1})$ for all $n\geq 1$. It follows that $\Lip(f)=2$. Now, a simple application of the Hypercyclicity Criterion for Lipschitz operators \cite[Theorem 3.1]{ACP} gives us that $\widehat{f}$ is hypercyclic. So, $R_{\widehat{f}}$ is dense in $\F(M)$.
\end{proof}
Further results about the strong operator topology and weak operator topology restricted to the set $\{\widehat{f}:~f\in \Lip_0(M,N)\}$ are given in Section~\ref{sec: top}. \\

The rest of this section is dedicated to finish the proof of Theorem~\ref{theo: main 2}, that is, to show that $\textup{int}(R_{\widehat{f}})\neq \emptyset$ implies that $R_{\widehat{f}}=\F(M)$ for the following two classes of metric spaces: complete metric spaces with at most countable many non-isolated points and (non necessarily bounded) connected subset of $\R$.
In order to do so, we need some partial results which may have their own interest. 
\begin{proposition}\label{recurrent no zero}
	Let $(M,d)$ be a pointed metric space and $f\in \Lip_0(M,M)$ such that $\textup{int}(R_{\widehat{f}})\neq \emptyset$.
	Then, there is no $\mu\in \mathcal{F}(M)\setminus \textup{span}(\delta(M))$ such that the sequence $(\widehat{f}^n(\mu))_n$ converges in $\textup{span}(\delta(M))$. 
	In particular, $\widehat{f}$ is injective. 
\end{proposition}
In order to prove Proposition~\ref{recurrent no zero}, we need the following lemma.
\begin{lemma}\label{closed sum}
Let $k\in \N$ and $(\lambda_i)_{i=1}^k\subset\R$. 
Then, $\sum_{i=1}^k \lambda_i \delta(M)$ is a norm closed subset of $\F(M)$.
\end{lemma}
\begin{proof}
Let us write $E:=\sum_{i=1}^k \lambda_i \delta(M)$.
Let $(u_n)_n\subset E$ and $\mu\in \F(M)$ such that $(u_n)_n$ converges to $\mu$ in norm. 
In order to fix notation, for $n\in\N$, let $(x_{i,n})_{i=1}^k\subset M$ such that $u_n=\sum_{i=1}^{k}\lambda_i \delta(x_{i,n})$.\\

Thanks to \cite[Lemma 2.10]{ACP}, $\textup{supp}(\mu)$ has at most $k$ elements. 
So, there is $N\subset M$, $|N|\leq k$ and $(\alpha_x)_{x\in N}\subset \R$ such that $\mu=\sum_{x\in N}\alpha_x\delta(x)$.
Since $N$ is a finite set, there are an increasing sequence $(n(j))_j\subset \N$, $\varepsilon>0$ and $I\subset\{1,...,k\}$ such that for each $i\in I$, the sequence $(x_{i,n(j)})_j$ converges to a point in $N$, and for each $i\notin I$,  $\textup{dist}(x_{i,n(j)},N)>\varepsilon$ for all $j\in\N$.
Now, consider $I_x:=\{i\in I:~\lim_jx_{i,n(j)}=x\}$ for any $x\in N$. 
Then, $\alpha_x=\sum_{i\in I_x}\lambda_i$ for all $x\in N$. 
Also, thanks to Proposition~\ref{prop: supp mu}, we get that $\lim_j \sum_{i\notin I} \lambda_i\delta(x_{i,n(j)})=0$.
This ends the proof of the Lemma~\ref{closed sum}.

\end{proof}

\begin{proof}[Proof of Proposition~\ref{recurrent no zero}]
	
	Let us assume that there is $\mu\in \mathcal{F}(M)\setminus \textup{span}(\delta(M))$ such that $(\widehat{f}^n(\mu))_n$ converges in norm to a vector in $\textup{span}(\delta(M))$.
	Let $u=\sum_{i=1}^k\lambda_i \delta(x_i)\in \textup{int}(R_{\widehat{f}})$, where $(\lambda_i)_{i=1}^k\subset\R$ and $(x_i)_{i=1}^k\subset M$.
	Thanks to Lemma~\ref{closed sum}, any accumulation point of the sequence $(\widehat{f}^n(u+\varepsilon\mu))_n$ belongs to $\sum_{i=1}^k \lambda_i \delta(M) + \textup{span}(\delta(M))= \textup{span}(\delta(M))$, and this for any $\varepsilon>0$. 
	Therefore, $u+\varepsilon\mu \notin R_{\widehat{f}}$ for any $\varepsilon>0$.
	Thus, $u\notin \textup{int} (R_{\widehat{f}})$, which is a contradiction.
	Let us see now that $\widehat{f}$ is injective. 
	Reasoning again towards a contradiction, let us assume that there is $u\in \textup{span}(\delta(M))$ such that $\widehat{f}(u)=0$. Therefore, $u\notin R_{\widehat{f}}$, which contradicts Corollary~\ref{cor: Rf dense}.
	So, there is $\mu\in \F(M)\setminus\textup{span}(\delta(M))$ such that $\widehat{f}(\mu)=0$. This contradicts the first part of this proposition. 
	Hence, $\widehat{f}$ is injective. 	
\end{proof}

Now, we turn our study to the case whenever $(M,d)$ is a complete metric space with at most countably many non-isolated points.

\begin{lemma}\label{lem: countable spaces}
	Let $(M,d)$ be a complete pointed metric space.
	Assume that $(M,d)$ has at most countably many non-isolated points.
	Let $f\in \Lip_0(M,M)$ such that $R_f=M$. 
	Then, for any $x\in M$ the set $\textup{orb}(f,x)$ is finite, i.e. every point in $M$ is a periodic point for $f$. 
\end{lemma}

\begin{proof}  
	Let us assume first that $M$ is countable.
	Let us write $M_1=M$.
	Since $M$ is countable and complete, thanks to the Baire category theorem we can find an isolated point $x\in M$. 
	Let $\varepsilon>0$ such that $B(x,\varepsilon)=\{x\}$.
	Since $x$ is a recurrent point for $f$, $x$ is a periodic point of $f$. 
	Also, the continuity of $f$ and the fact that $R_f=M$ imply that $\textup{orb}(f,x)$ contains only isolated points of $M$. 
	Indeed, otherwise there are a non-isolated point $y\in \textup{orb} (f,x)$ and $n\in \N$ such that $f^n(y)=x$. 
	Thus, by continuity of $f$ there is $\rho'>0$ such that $f^n(B(y,\rho'))\subset B(x,\varepsilon)=\{x\}$, where $B(y,\rho')$ is an infinite set. 
	This contradicts the fact that $R_f=M$.
	Thus, setting $I_1=\textup{orb}(f,x)$, we have that every point in $I_1$ is isolated, $M_2:=M\setminus I_1$ is a closed subset of $M$ and $f(M_2)\subset M_2$. 
	Indeed, the last inclusion follows from the fact that $R_f=M$.\\

	Let us proceed by transfinite induction. Let $\alpha$ be an ordinal and assume that we have constructed a collection of finite sets $\{I_\beta:~\beta<\alpha\}$ satisfying: for any $\beta<\alpha$ and any $x\in I_\beta$
	
	\begin{enumerate}
		\item[$i)$] $M_\beta:= M\setminus\bigcup\{I_\gamma:~\gamma<\beta\}$ is closed,
		\item[$ii)$] $f(M_\beta)\subset M_\beta$, 
		\item[$iii)$] $x$ is an isolated point of $M\setminus\bigcup\{I_\gamma:~\gamma<\beta\}$, and
		\item[$iv)$] $I_\beta=\textup{orb}(f,x)$.
	\end{enumerate}
	Now, note that the set
	\[M_\alpha:=\bigcap_{\beta<\alpha} M_\beta\]
	is closed for being an arbitrary intersection of closed sets.
	If $M_\alpha$ is empty, there is nothing to do. 
	If $M_\alpha\neq\emptyset$, thanks to the statement $ii)$ above,  $f(M_\alpha)\subset M_\alpha$. Thus, thanks to the Baire category theorem, we get that there is a point $x_\alpha\in M_\alpha$ which is isolated in $M_\alpha$. 
	Following the same argument as above, we get that $x_\alpha$ is a periodic point for $f$ and $I_\alpha:=\textup{orb}(f,x)$ only contains isolated points of $M_\alpha$. This complete the induction step.\\
	
	Observe that, since $I_\alpha$ is nonempty if $M_\alpha$ is nonempty, the induction ends at some countable ordinal $\Lambda$.
	It follows that $M=\bigcup\{I_\beta:~\beta<\Lambda\}$ and thus, all the orbits generated by the action of $f$ are finite.\\
	
	Now, if $M$ is uncountable, let us write by $N$ the set of isolated points of $M$. Reasoning as the previous case, we get that $f(N)\subset N$, $f(M\setminus N)\subset M\setminus N$ and that $N$ only contains periodic points for $f$. Now, noticing that $M\setminus N$ is a countable, closed subset of $M$ and $f$-invariant, the first part of this proof show that $M\setminus N$ also contains only periodic points for $f$.
\end{proof}

Albeit simple, the following lemma is crucial for our result. The details are left to the reader.

\begin{lemma}\label{lem: acc of acc}
	Let $(M,d)$ be a complete metric space and let $f:M\to M$ be a continuous function. Let $x\in M$ and $y$ be an accumulation point of $(f^n(x))_n$. Then, every accumulation point of $(f^n(y))_n$ is an accumulation point of $(f^n(x))_n$.
\end{lemma}

Now we can prove Theorem~\ref{theo: main 2} (1).

\begin{theorem}\label{countable spaces}	
	Let $M$ be a pointed metric space and $f\in \Lip_0(M,M)$ such that  $\textup{int}(R_{\widehat{f}})\neq \emptyset$.
	Assume further that every orbit of $f$ is finite. 
	Then, $R_{\widehat{f}}=\F(M)$.
\end{theorem}

Followed by Lemma~\ref{lem: countable spaces} and Theorem~\ref{countable spaces}, we have the following consequence.
\begin{corollary}
	Let $(M,d)$ be a complete pointed metric space with at most countable many non-isolated points.
	Then, for any $f\in \Lip_0(M,M)$ the following dichotomy holds true: either $R_{\widehat{f}}$ has empty interior or $R_{\widehat{f}}=\F(M)$.  
\end{corollary}
 
\begin{proof}[Proof of Theorem~\ref{countable spaces}]
	Let $f\in \Lip_0(M,M)$ such that $\textup{int}(R_{\widehat{f}})\neq \emptyset$ and $\textup{orb}(f,x)$ is finite for every $x\in M$. 
	Thanks to Proposition~\ref{prop: recurrent points interior}, $R_f=M$. Thus, every point in $M$ is a periodic point for $f$. 
	Let $u\in \textup{span}(\delta(M))$ and $\varepsilon>0$ such that $B(u,\varepsilon)\subset \F(M)$. 
	Thus, there are an one-to-one finite sequence $(x_i)_{i\in I}\subset M\setminus\{0\}$ and $(\lambda_i)_{i\in I}\subset \R\setminus\{0\}$ such that $u=\sum_{i\in I } \lambda_i\delta(x_i)$.
	Let $M_u=\textup{supp}(u)$.
	Since every point has a finite orbit, by shrinking $\varepsilon$ and perturbing $u$ if necessary, we assume that 
	\[M_u=\bigcup_{x\in M_u}\textup{orb}(f,x).\]
	Observe that $f(M_u)= M_u$ and that there is $\sigma:I\to I$ bijective such that $f(x_i)=x_{\sigma(i)}$ for all $i\in I$. 
	\\
	
	We prove that $B(0,\varepsilon)\subset R_{\widehat{f}}$.
	To this end, let $\mu\in \F(M)$ with $\|\mu\|< \varepsilon$. 
	Thus, $u+\mu\in R_{\widehat{f}}$. 
	Therefore, there is $(n(j))_j\subset \N$, such that $\widehat{f}^{n(j)}(u+\mu)$ tends to $u+\mu$ as $j$ tends to infinity. 
	Since $M_I$ is a finite set and $f$ acts like a permutation on $M_I$, up to subsequence of $(n(j))_j$ which will be still denoted by $(n(j))_j$, we have that $f^{n(j)}(x_i)=x_{\pi(i)}$ for all $i\in I$ and all $j\in \N$, where $\pi:I\to I$ is a bijective function.
	Let us write $u_{\pi}:= \sum_{i\in I}\lambda_ix_{\pi(i)}$. 
	With this, we have that 
	
	\[\lim_{j\to\infty} \widehat{f}^{n(j)}(\mu)=\lim_{j\to\infty} \widehat{f}^{n(j)}(\mu+u)-\widehat{f}^{n(j)}(u)= \mu+u-u_\pi. \]

	Note that, for any $m\geq 1$, we have that
	\[\lim_{j\to\infty} \widehat{f}(\mu+u-u_{\pi^m})=(\mu+u-u_\pi) +u_\pi -u_{\pi^{m+1}}=\mu+u-u_{\pi^{m+1}}.\] 	
	
	Now, since $I$ is a finite set, there is $\alpha\in \N$ such that $\pi^{\alpha}$ is the identity on $I$ (for instance, $\alpha= |I|!$). Therefore, $\mu+u-u_{\pi^{\alpha}}=\mu$. 
	Finally, applying repeatedly Lemma~\ref{lem: acc of acc}, we get that $\mu$ is an accumulation point of $(\widehat{f}^j(\mu))_j$ and thus, $\mu\in R_{\widehat{f}}$. 
	Since $R_{\widehat{f}}$ is invariant under non-zero scalar multiplication, we get that $R_{\widehat{f}}=\F(M)$. 
\end{proof}

Let us end this section with the case whenever the underlying metric spaces is a closed connected subset of $\R$ (Theorem~\ref{theo: main 2}, $ii$)). 
The case of compact intervals can be also obtained as a consequence of Proposition~\ref{recurrent no zero} and the nice result of Block and Coven given in \cite[Theorem]{BC}. Here we present an independent proof which can be applied to the unbounded case.
\begin{theorem}
	Let $I\subset \R$ be a (not necessarily bounded) closed connected set such that $0\in I$.
	Let $f\in\Lip_0(I,I)$ such that $\textup{int}(R_{\widehat{f}})\neq \emptyset$. 
	Then, $f^2$ is the identity on $I$. In particular, $R_{\widehat{f}}=\F(I)$.
\end{theorem}
\begin{proof}
	Since $R_{\widehat{f}}$ has nonempty interior, by Proposition~\ref{recurrent no zero}, we know that $\widehat{f}$ is injective. 
	Thus, $f$ is injective as well and therefore, either increasing or decreasing. \\
	
	Case 1: Let us assume that $f$ is an increasing function. We show that $f=id$. Indeed, otherwise there is $x\in I$ such that $f(x)\neq x$. 
	Let us assume that $f(x)>x$. 
	Since $f$ is increasing, $f^n(x)\geq f(x)>x$ for all $n\in \N$. Thus, $x\notin R_f$ implying that $\delta(x)\notin R_{\widehat{f}}$. 
	This contradicts Proposition~\ref{prop: recurrent points interior}.
	On the other hand, if $f(x)<x$, since $f$ is increasing, we have that $f^n(x)\leq f(x)<x$ for all $n\in \N$. We get the same contradiction. Hence, $f(x)=x$ for all $x\in \R$.\\
	
	Case 2: Let us assume that $f$ is a decreasing function. Then, $f((0,\infty)\cap I)=\subset (-\infty,0)$ and $f((-\infty,0)\cap I)\subset(0,\infty)$. 
	Observe that, for any $x\in \R$, if $(f^{n(j)}(x))_j$ converges to $x$, there is $J\in \N$ such that $n(j)$ is even for all $j\geq J$. 
	With this and Proposition~\ref{prop: recurrent points interior}, we have that $f^2$ is an increasing function such that $R_{f^2}=I$.
	We can apply the same argument of Case 1 to conclude that $f^2=id$.\\
	
	Finally, we have that $\widehat{f}^{2n}=Id$ on $\F(I)$ for all $n\in\N$.
	Thus, $R_{\widehat{f}}=\F(I)$. 
\end{proof}

\section{Relation between $A_{\widehat{f}}$ and $U_f$: two examples}\label{sec: U set}
Let $(M,d)$ be a pointed metric space and $f\in \Lip_0(M,M)$.
The aim of this section is to understand the relation between the sets $A_f$, $U_f$, $A_{\widehat{f}}$ and $U_{\widehat{f}}$. 
In fact, this section is motivated by the following question: Is there a canonical way to construct a vector $\mu\in A_{\widehat{f}}$?.
The following diagram contains all the implications which (trivially) hold true. 
In what follows, we show that all the unwritten implications are false in general.

\begin{center}
\includegraphics[scale=0.45]{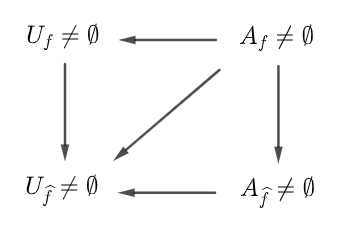}	
\end{center}

Our first example, Corollary~\ref{cor: comparison Uf Af}, shows that there is a metric space $(M,d)$ and a function $f\in \Lip_0(M,M)$ such that $U_{\widehat{f}}=\F(M)\setminus \{0\}$ but $A_{\widehat{f}}=\emptyset$.
Then, in Example~\ref{example 2}, we explicit a bounded metric space $(M,d)$ and a function $f\in \Lip_0(M,M)$ such that $A_{\widehat{f}}\neq \emptyset$. Since $M$ is bounded, it follows that $U_f$ and $A_f$ are necessarily empty.\\  

Let us start with the preliminaries for the first example.
\begin{proposition}\label{prop: alpha}
	There are an unbounded sequence $(\alpha_n)_n\subset \R_+$ and $L>0$ such that $\alpha_{n+1}\leq L\alpha_n$ for all $n\geq 1$ and 
 	\begin{align}\label{alpha n property}
		\liminf_{m\to\infty} \sum_{n=1}^\infty |\lambda_n| \dfrac{\alpha_{n+m}}{\alpha_n}<\infty,~ \forall \lambda\in \ell^1.
	\end{align}
\end{proposition}

\begin{proof}
	Let $(k_n)_n\subset\mathds{N}$ be an increasing sequence with $k_0=0$. 
	Let us define $s_n=\sum_{i=0}^n k_i$ for any $n\geq 0$. 
	Assume that, for any $m\geq 1$, we have 
	\[m^{\frac{s_{m-1}}{k_m}}\leq 2,~~~\hspace{0.8cm}~\text{and}\hspace{0.5cm} m^{1/k_m}\geq (m+1)^{1/k_{m+1}}.\]
	Now, for any $m\in\N$ and $n\in (s_m,s_{m+1}]$, let us write $n:=s_m+j$ and define 
	\[\alpha_n=\alpha_{s_m+j}=(m+1)^{j/k_{m+1}}.\]
	Observe that $\alpha_{s_m}=m$ for all $m\geq 1$.
	We claim that the sequence $(\alpha_n)_n$ witnesses the inequality~\eqref{alpha n property}.
	Indeed, it is clear that $(\alpha_n)_n$ is an unbounded sequence. 
	Now, let us fix $m\in \mathds{N}$.
	For $n\in (s_m,\infty)$ we have that $\alpha_{n+s_m}\leq  (m+1)^{\frac{s_m}{k_{m+1}}}\alpha_n\leq 2\alpha_n$. 
	On the other hand, if $n\leq s_m$, then $s_m\leq n+s_m\leq 2s_m$ and thus $\alpha_{n+s_m}\leq 2$. Since $\alpha_n\geq 1$ for all $n$, we get that $\alpha_{n+s_m}/\alpha_n\leq 2$ for all $m\in \N$ and all $n\in \N$.
	Summarizing, we have that for any $\lambda:=(\lambda_n)_n\in \ell^1$

	\[\liminf_{m\to\infty} \sum_{n=1}^\infty |\lambda_n| \dfrac{\alpha_{n+m}}{\alpha_n}\leq \liminf_{m\to\infty} \sum_{n=1}^\infty |\lambda_n| \dfrac{\alpha_{n+s_m}}{\alpha_n}\leq 2\| (\lambda_n)_n\|_1.\]
	Finally, considering $L=2$, the proof is complete. 
\end{proof}
\begin{remark}
	In the proof of Proposition~\ref{prop: alpha}, a stronger property is proved. 
	Indeed, we found sequences $(\alpha_n)_n\subset\R_+$ and $(m_k)_k\subset \N$ such that $(\alpha_n)_n$ is unbounded and the set $\{(\alpha_{n+m_k}/\alpha_n)_n:~k\in\N\}$ is bounded in $\ell^\infty$. This clearly implies~\eqref{alpha n property}.
\end{remark}

Let $\alpha=(\alpha_n)_n\subset \R_+$ and $f:(\N,d_\alpha)\to(\N,d_\alpha)$ defined by $f(0)=0$ and $f(n)=n+1$ for all $n\geq 1$.
Thanks to Proposition~\ref{prop: function on N alpha}, the function $f$ is $L$-Lipschitz if and only if $\alpha_{n+1}\leq L\alpha_n$ for all $n\geq 1$.  
In what follows we study the orbits of the function $f$ with respect different choices of $\alpha$. 
In \cite[Remark 1.3]{ACP} it is shown that if $\alpha$ is bounded away from $0$ and infinity, then $\widehat{f}$ only has bounded orbits. Indeed, let $r,R>0$ such that $r<\alpha_n<R$ for all $n\in \N$. Let $\mu\in \mathcal{F}(\N)$. Thanks to Proposition~\ref{free of N}, there is $\lambda=(\lambda_n)_n\in \ell^1(\N)$ such that $\mu=\sum_{n=1}^\infty \lambda_n\delta(n)/\alpha_n$.
Thus, we have that
\[\| \widehat{f}^j(\mu)\| = \sum_{n=1}^\infty |\lambda_n|\dfrac{\|\delta(f^j(n))\|}{\alpha_n}\leq \dfrac{R}{r}\|\lambda\|_1,~\forall j\in\N.\]

\begin{example}\label{example 1}
	Let $\alpha=(\alpha_n)_n\subset \R_+$ and $f:(\N,d_\alpha)\to(\N,d_\alpha)$ defined by $f(0)=0$ and $f(n)=n+1$ for all $n\geq 1$. Further, assume that $f$ is a Lipschitz function.
	\begin{enumerate}
		
		\item[a)] If $\alpha$ satisfies \eqref{alpha n property}, then $A_{\widehat{f}}=\emptyset$. 
		\item[b)] If $\alpha$ do not satisfy \eqref{alpha n property}, then $A_{\widehat{f}}$ is dense in $\mathcal{F}(\N)$.
	\end{enumerate}
	Moreover, if $\alpha$ is unbounded, $U_{\widehat{f}}=\mathcal{F}(\N)\setminus\{0\}$.
\end{example}

\begin{proof}
	\textbf{a):} Let $\mu\in \mathcal{F}(\N)$. 
	Then, thanks to Proposition~\ref{free of N}, there is $\lambda=(\lambda_n)_n\in \ell^1(\N)$ such that $\mu=\sum_{n=1}^\infty \lambda_n\delta(n)/\alpha_n$. 
	Now, thanks to \eqref{alpha n property}, we have that
	\[\liminf_{k\to\infty} \|\widehat{f}^k(\mu)\| = \liminf_{k\to\infty} \sum_{n=1}^\infty |\lambda_n|\dfrac{\alpha_{k+n}}{\alpha_n}<\infty.\]
	
	So, $\mu\notin A_{\widehat{f}}$. \\
	\textbf{b):}
	Let $\lambda\in \ell^1(\N)$ such that the inequality~\eqref{alpha n property} is not satisfied. 
	Let $\mu=\sum_{n=1}^\infty \lambda_n\frac{\delta(n)}{\alpha_n}\in \mathcal{F}(\N)$. Let us see that $\mu\in A_{\widehat{f}}$. 
	Indeed,
	\[\liminf_{k\to\infty} \|\widehat{f}^k(\mu)\|= \liminf_{k\to\infty} \sum_{n=1}^\infty |\lambda_n|\dfrac{\alpha_{k+n}}{\alpha_n}=\infty. \]
	Now, we prove that $\textup{span}(\delta(\N))\subset \overline{A}_{\widehat{f}}$.
	To this end, let $u= \sum_{n=1}^\infty \beta_n\delta(n)\in \textup{span}(\delta(M))$, i.e. $(\beta_n)_n$ is a finitely supported sequence. 
	Let $N>0$ such that $\beta_n=0$ for all $n\geq N$.
	Let $\varepsilon >0$ such that $|\beta_n+\varepsilon\lambda_n|\geq \varepsilon|\lambda_n|$ for all $n\in \{1,...,N\}$. 
	Observe that $\varepsilon$ can be chosen arbitrarily close to $0$. 
	Thus
	\begin{align*}
	\liminf_{k\to\infty} \|\widehat{f}^k(u+\varepsilon \mu)\|&=\liminf_{k\to\infty} \sum_{n=1}^N |\beta_n+\varepsilon \lambda_n|\dfrac{\alpha_{k+n}}{\alpha_n}+\sum_{n=N+1}^\infty \varepsilon|\lambda_n|\dfrac{\alpha_{k+n}}{\alpha_n}\\
	&\geq \varepsilon \liminf_{k\to\infty} \sum_{n=1}^\infty |\lambda_n|\dfrac{\alpha_{k+n}}{\alpha_n}=\infty.
	\end{align*}
	
	Thus $u+\varepsilon\mu \in A_{\widehat{f}}$.
	Now, by sending $\varepsilon$ to $0$, we get that $u\in \overline{A}_{\widehat{f}}$ as claimed. 
	Finally, we observe that $\mathcal{F}(M)=\overline{\textup{span}}(\delta(\N))\subseteq \overline{A}_{\widehat{f}}$.\\

	It only remains to prove the last assertion of the proposition. 
	Let $\alpha$ be an unbounded sequence. 
	Let $\mu\in \mathcal{F}(\N)$ different from $0$ and $\lambda\in \ell^1(\N)$ such that $\mu = \sum_{n=1}^\infty\lambda_n\delta(n)/\alpha_n$. 
	Let $m\in \N$ such that $\lambda_m\neq 0$. 
	Then, 
	\[\limsup_{k\to\infty} \|\widehat{f}^k(\mu)\| = \limsup_{k\to\infty}\sum_{n=1}^\infty |\lambda_n|\dfrac{\|\delta(f^k(n))\|}{\alpha_n}\geq\dfrac{|\lambda_m|}{\alpha_m} \limsup_{k\to\infty} \alpha_{m+k}=\infty.\]
	So, $\mu\in U_{\widehat{f}}$. 
\end{proof}

\begin{corollary}\label{cor: comparison Uf Af}
	There is a pointed metric space $(M,d)$ and a function $f\in \textup{Lip}_0(M,M)$ such that $U_{\widehat{f}}=\mathcal{F}(M)\setminus\{0\}$ but $A_{\widehat{f}}=\emptyset$.
\end{corollary}
\begin{proof}
	It is a direct consequence of Proposition~\ref{prop: alpha} and Example~\ref{example 1} $a)$.
\end{proof}

Let us now present our second example.
\begin{example}\label{example 2}
	Let $f:[0,1]\to [0,1]$ defined by $f(x)=\min(2x,1)$. Then, $U_f=\emptyset$ and $A_{\widehat{f}}$ is nonempty.
\end{example}
\begin{proof}
	Since $[0,1]$ is a bounded metric space, $U_f=A_f=\emptyset$. 
	Observe also that $\textup{Lip}(f)=2$. 
	For $n\in\N$, let $x_n=2^{-n}$.  Since $\|\delta(x_n)\|=2^{-n}$ for all $n$, we have that $\mu:= \sum_{n=1}^\infty \delta(x_n)
	\in \mathcal{F}([0,1])$.
	Observe that, for any $k\in \N$, $f^k(\mu)= k\delta(1)+\mu$. 
	So, $\mu\in A_{\widehat{f}}$.
\end{proof}
\begin{remark}
	In fact, in the context of Example~\ref{example 2}, $A_{\widehat{f}}$ is dense in $\mathcal{F}([0,1])$ (see Proposition~\ref{bounded orbits}).
\end{remark}

\section{Vectors escaping to infinity}\label{sec: A set}

In this section we focus our study on the set of vectors escaping to infinity for a given Lipschitz operator, i.e. $A_{\widehat{f}}$.
In fact, we deal with the following question:
is it true that $A_{\widehat{f}}$ is always either empty or dense in $\mathcal{F}(M)$? 
Recall that $A_{\widehat{f}}:=\{\mu\in \F(M):~\lim_n\|\widehat{f}^n(\mu)\|=\infty\}$.\\

In what follows, we collect several results pointing towards a positive answer of the mentioned question.
We will provide a complete proof in the case where $M$ is a closed connected subset of $\R$. 
The question is still open in its full generality.\\

Let us start with the following result.

\begin{proposition}\label{bounded orbits}
	Let $M$ be a pointed metric space. Let $f\in \Lip_0(M,M)$ such that every $\textup{orb}(f,x)$ is bounded for all $x\in M$. Then, $A_{\widehat{f}}$ is empty or dense in $\mathcal{F}(M)$.
\end{proposition}

\begin{proof}
	Let us assume that there is $\mu\in A_{\widehat{f}}$. 
	It is enough to prove that $\mu+\textup{span}(\delta(M))\subset A_{\widehat{f}}$.\\
	
	Let $u\in \textup{span}(\delta(M))$. 
	There are $(x_i)_{i=1}^{k}\subset M$ and $(\lambda_i)_{i=1}^{k}\subset\R$ such that
	$u=\sum_{i=1}^{k} \lambda_i \delta(x_i)$.
	Let $0<C<\infty$ such that 
	$\textup{diam}(\{0\}\cup \textup{orb}(f,x_i))\leq C$ for all $i=1,...,k$.
	Therefore, for any $n\in \N$, we have that
	\[ \| \widehat{f}^n (u+\mu)\|\geq \|\widehat{f}^n(\mu)\| -  \sum_{i=1}^{k} |\lambda_i| \|\delta(f^{n}(x_i))\| = \|\widehat{f}^n(\mu)\| -  C\sum_{i=1}^{k} |\lambda_i|. \]
	
	Since $\mu\in A_{\widehat{f}}$, $\lim_{n\to\infty}\| \widehat{f}^n (u+\mu)\|=\infty$. This proves that $A_{\widehat{f}}\supset \mu+\textup{span}(\delta(M))$ and the proof is complete.
\end{proof}

In the following proposition we state another sufficient condition to obtain the density of the set $A_{\widehat{f}}$.

\begin{proposition}\label{prop: gen bounded space}
	Let $M$ be a pointed metric space. 
	Let $f\in \Lip_0(M,M)$ and $\mu\in A_{\widehat{f}}$. Assume that there is a sequence $(\varphi_n)_n\in \Lip_0(M)$ of $1$-Lipschitz functions such that $\bigcup_n\textup{supp}(\varphi_n)$ is a bounded subset of $M$ and $(\langle \varphi_n,f^n(\mu)\rangle)_n$ tends to infinity. Then, $A_{\widehat{f}}$ is dense in $\mathcal{F}(M)$.
\end{proposition}
\begin{proof}
	We show that $\mu+\textup{span}(\delta(M))\subset A_{\widehat{f}}$. 
	Let $R>0$ such that $\bigcup_n \textup{supp}(\varphi_n)\subset B(0,R)$. 
	Let $u\in \textup{span}(\delta(M))$. 
	Thus, there are $(x_i)_{i=1}^k\subset M$ and $(\lambda_i)_{i=1}^k\subset \R$ such that $u=\sum_{i=1}^k\lambda_i\delta(x_i)$.
	Now, note that 
	\[\|\widehat{f}^n(\mu+u)\| \geq \langle \varphi_n, \widehat{f}^n(\mu +u)\rangle\geq \langle \varphi_n,\widehat{f}^n(\mu)\rangle - R\sum_{i=1}^k|\lambda_i|,\]
	which tends to infinity as $n$ tends to infinity. Thus, $u+\mu\in A_{\widehat{f}}$.
\end{proof}

Proposition~\ref{prop: gen bounded space} and Hahn-Banach Theorem readily imply the following result.
\begin{corollary}\label{cor: boundedness of fn(mu)}
	Let $M$ be a complete pointed metric space. 
	Let $f\in\Lip_0(M,M)$ and $\mu\in A_{\widehat{f}}$. Assume that $\bigcup_n \textup{supp}(\widehat{f}^n(\mu))$ is bounded. Then $A_{\widehat{f}}$ is dense in $\F(M)$.
\end{corollary}

Let us continue with the main (abstract) result of this section.

\begin{theorem}\label{theo: int af implies af dense}
	Let $M$ be a pointed metric space. 
	Let $f\in \Lip_0(M,M)$ and assume that there is $v\in \textup{span}(\delta(M))\cap  A_{\widehat{f}}$. 
	Then, $A_{\widehat{f}}$ is dense in $\mathcal{F}(M)$.
\end{theorem}

\begin{proof}
	Let $(x_i)_{i=1}^k\subset M\setminus\{0\}$ be a one-to-one sequence and $(\lambda_i)_{i=1}^k\subset \R\setminus\{0\}$ such that $v=\sum_{i=1}^k\lambda_i\delta(x_i)$.
	Since $v\in A_{\widehat{f}}$, without loss of generality we can assume that $\textup{orb}(f,x_i)$ is unbounded for all $i\in \{1,...,k\}$. Moreover, we can (and shall) assume that $\lambda_i=1$ for all $i=1,...,k$. Indeed, since $v=\sum_{i=1}^k\lambda_i\delta(x_i)\in A_{\widehat{f}}$, for any $R>0$ there is $N\in \N$ such that $\|\widehat{f}^n(v)\| \geq R\sum_{i=1}^k|\lambda_i|$ for all $n\geq N$. 
	Thus, there is at least one $i\in\{1,...,k\}$ (depending on $n$) such that $d(f^n(x_i),0)\geq R$ for all $n\geq N$. 
	Therefore $\|\widehat{f}^n(\sum_{i=1}^k \delta(x_i))\|\geq R$ for all $n\geq N$. From now on, $\lambda_i=1$ for all $i\in\{1,...,k\}$.
	\\

	We prove that $\textup{span}(\delta(M))\subset \overline{A}_{\widehat{f}}$. 
	To this end, let $u=\sum_{i=1}^l\alpha_i\delta(y_i)\in \textup{span}(\delta(M))$, where $\alpha_i\in \R\setminus\{0\}$ and $y_i\in M$ for all $i=1,...,l$. 
	Let us consider $\{I,~J\}$ be a partition of $\{1,...,l\}$ such that $\textup{orb}(f,y_i)$ is bounded for all $i\in I$ and unbounded for all $i\in J$. 
	Let $R>0$ such that $R>\textup{diam}(\{0\}\cup\bigcup_{i\in I}\textup{orb}(f,y_i))$.
	Also, let $\varepsilon>0$ and $w=u+\varepsilon v$.
	Since $v\in A_{\widehat{f}}$, there is $N\in \N$ such that for all $n\geq N$, there is $i\in \{1,...,k\}$ such that $d(f^n(x_i),0)\geq (|J|+k+2)R$.
	Let us fix $n\geq N$.
	Note that the cardinal of $\{x_i,y_j:~i=1,...,k,~j\in J\}$ is at most $|J|+k$ and at least one of them that does not belong to $B(0, (|J|+k+2)R)$. 
	By the pigeonhole principle,
	there is $m\in \{1,...,|J|+k+1\}$ such that
	\[(\{f^n(x_i):~i\in \{1,...,k\}\}\cup\{f^n(y_j):~j\in J\})\cap B(0,(m+1)R)\setminus \overline{B}(0,mR)=\emptyset. \]
	
	Let $\varphi\in \textup{Lip}_0(M)$ be the function defined by
	\[\varphi(x)= \min (\max (d(x,0)-mR,0),R).\]
	Observe that $\Lip(\varphi)=1$ and $0\leq \varphi\leq R$. 
	Moreover, $\varphi(f^n(x_i))\in \{0,R\}$ for all $i=1,...,k$, $\varphi(f^n(y_i))=0$ if $j\in I$, and $\varphi(f^n(y_i))\in\{0,R\}$ if $i\in J$.
	Now, note that
	\begin{align*}
		\|\widehat{f}^n(v)\|&\geq |\langle \varphi,\widehat{f}^n(u)+\varepsilon\widehat{f}^n(\nu)\rangle| \\
		&=\left| \left\langle \varphi,\sum_{i\in J} \alpha_i\delta(f^n(y_i))+\varepsilon\sum_{i=1}^k \delta(f^n(x_i))\right\rangle\right|\geq RC_\varepsilon,
	\end{align*}
	where $C_\varepsilon$ can be chosen as a constant which only depends on $\varepsilon$ and $u$, for instance:
	\[C_\varepsilon := \min \left\{\left |\sum_{i\in L}\alpha_i+ \varepsilon\rho \right |:~L\subset J,~\rho\in \{1,...,k\} \right\}.\]
	Since $J$ is a finite set, there are arbitrarily small $\varepsilon>0$ such that $C_\varepsilon>0$. 
	Let us fix $\varepsilon>0$ such that $C_\varepsilon>0$. 
	So, we have that $\|\widehat{f}^n(v)\|\geq RC_\varepsilon$ for any $n\geq N$. 
	Finally, since $R$ can be chosen arbitrarily large, we get that $u+\varepsilon v\in A_{\widehat{f}}$. 
	Also, since $C_\varepsilon$ is greater than $0$ for arbitrarily small $\varepsilon$, we get that $u\in \overline{A_{\widehat{f}}}$. 
\end{proof}

As a direct consequence of Theorem~\ref{theo: int af implies af dense} we have:
\begin{corollary}\label{cor: examples 2}
Let $M$ be a pointed metric space and let $f\in \Lip_0(M,M)$ such that
\begin{enumerate}
	\item[$i)$] $\textup{int}(A_{\widehat{f}})\neq \emptyset $, or
	\item[$ii)$] there is $x\in A_f$.

\end{enumerate}
Then, $A_{\widehat{f}}$ is dense. 
In particular, there is no wild operators of the form $\widehat{f}$.
\end{corollary}
The following corollary gather some examples which are consequence of Proposition~\ref{bounded orbits} and Theorem~\ref{theo: int af implies af dense}.

\begin{corollary}\label{cor: examples}
	The following pointed metric spaces $M$ satisfy the following dichotomy: for any $f\in \Lip_0(M,M)$, $A_{\widehat{f}}$ is either empty or dense in $\mathcal{F}(M)$.
	\begin{enumerate}
		\item[$i)$] $M$ is a bounded metric space.
		\item[$ii)$] Each connected component of $M$ is bounded and, for any $R>0$, $B(0,R)$ intersects only finitely many connected components of $M$.
		\item[$iii)$] $M$ is a subset of $\mathds{Z}^d$, with $d\geq 1$.		
	\end{enumerate}
\end{corollary}
\begin{proof}
	Indeed, for $i)$ we directly apply Proposition~\ref{bounded orbits}. For $ii)$, it is enough to note that for any $f\in \Lip_0(M,M)$ and any $x\in M$,
	the sequence $(f^n(x))_n$ is bounded (equivalently, it visits a finite number of connected components of $M$) or it tends to infinity (it visits an infinite number of connected components of $M$). 
	Thus, for any $x\in M$, $\textup{orb}(f,x)$ is bounded or $x\in A_f$. 
	Hence, the desired conclusion follows from Proposition~\ref{bounded orbits} or Corollary~\ref{cor: examples 2} $ii)$ respectively.
	Finally, $iii)$ is a particular case of $ii)$.
\end{proof}

In what follows, we study the case where $M$ is a closed and connected subset of $\R$, Theorem~\ref{theo: main 3} ($iii$).

\begin{theorem}\label{theo: Af for R}
	Let $M$ be a closed and connected subset of $\R$. 
	Then, for any $f\in \Lip_0(M,M)$, the set $A_{\widehat{f}}$ is either empty or dense in $\F(M)$.	
\end{theorem}

If $M$ is a compact interval, then Theorem~\ref{theo: Af for R} is a direct consequence of Corollary~\ref{cor: examples} $i)$.
Let us continue with the case where $M=[a,\infty)$, with $a\leq 0$. The case where $M=(-\infty,a]$, with $a\geq 0$ is completely analogous.

\begin{proposition}\label{prop: semi infinite interval}
Let $a\leq 0$ and $M=[a,\infty)$ equipped with its usual metric induced by $\R$. Then, for any $f\in \Lip_0(M,M)$, the set $A_{\widehat{f}}$ is either empty or dense in $\F(M)$. 
\end{proposition}

\begin{proof}
Let $f\in \Lip_0(M,M)$. 
Let $\mathcal{U}:=\{x>0:~ f(x)>x\}$.
We split our analysis in three cases.\\

Case 1: $\mathcal{U}=\emptyset$. 
Let $x\in M$. Note that the sequence $(f^n(x))_n$ is bounded.
Indeed, let $A= \max\{|f(y)|:~y\in [a,0]\}$.
Then, $(f^n(x))_n\subset [a, \max(A,x)]$. 
Thus, $f$ only generates bounded orbits. 
So, Proposition~\ref{bounded orbits} gives us the desired conclusion.\\

Case 2: $\mathcal{U}\neq\emptyset$ and $\mathcal{U}$ contains an unbounded connected component. 
Let $x>0$ such that $[x,\infty)\subset\mathcal{U}$.
We claim that $x\in A_{f}$. 
Indeed, since $x\in \mathcal{U}$, we have that $x< f(x)\in\mathcal{U} $. 
Moreover,  $(f^n(x))_n$ is an strictly increasing sequence.
Reasoning towards a contradiction, let us assume that $\overline{x}=\lim_n f(x)<\infty$.
However, since $f$ is continuous and $f(\overline{x})>\overline{x}$, there is $\varepsilon>0$ such that $f(y)>\overline{x}$ for all $y\in (\overline{x}-\varepsilon,\overline{x}+\varepsilon)$. This clearly contradicts the fact that $(f^n(x))_n$ converges to $\overline{x}$. 
So $x\in A_f$.
Corollary~\ref{cor: examples 2} $ii)$ gives us that $A_{\widehat{f}}$ is dense in $\F(M)$.\\

Case 3: $\mathcal{U}\neq\emptyset$ and $\mathcal{U}$ does not contain an unbounded connected component.
Since $\mathcal{U}$ is an open subset of $(0,\infty)$, every connected component of $\mathcal{U}$ is an open interval.
Let $0\leq b<c$ such that $(b,c)$ is a connected component of $\mathcal{U}$.
Since $f$ is continuous, we necessarily have that $f(b)=b$ and $f(c)=c$.
Let us split our analysis in two more cases:\\

Case 3.1: $f((b,c))=(b,c)$. Let us prove that, for any $x\in (b,c)$, $(f^n(x))_n$ converges to $c$.
Indeed, let $x\in (b,c)$. Since $f((b,c))=(b,c)$, we have that $(f^n(x))_n\subset(b,c)$. 
Also, since $(b,c)\subset \mathcal{U}$, the sequence $(f^n(x))_n$ is increasing, thus, convergent. 
If $\overline{x}=\lim_{n}f^n(x)<c$, we get a contradiction as we did in Case 2.
Now, let $(x_n)_n\subset (b,c)$ be a decreasing sequence such that $\sum_{n=1}^\infty |x_n-b|<\infty$.
Thus, $\mu:= \sum_{n=1}^\infty \delta(x_n)-\delta(b)$ absolutely converges and it belongs to $\F(M)$. 
Using the map $\varphi(\cdot)= \textup{dist}(\cdot,[a,b])\in \Lip_0(M)$, which satisfies $\Lip(\varphi)=1$, we get that
\[\|\widehat{f}^j(\mu)\|\geq \langle \varphi ,\widehat{f}^j(\mu)\rangle=\sum_{n=1}^\infty |f^j(x_n)-b|,~ \forall j\in \N. \]
Now, it easily follows that $(\|\widehat{f}^j(\mu)\|)_j$ tends to infinity as $j$ tends to infinity. 
So, $\mu\in A_{\widehat{f}}$. 
Also, we know that $\bigcup_n\textup{supp}(\widehat{f}^n(\mu))\subset [b,c]$.
Thus, thanks to Corollary~\ref{cor: boundedness of fn(mu)}, $A_{\widehat{f}}$ is dense in $\F(M)$.\\

Case 3.2: $f((b,c))\supsetneq(b,c)$. 
Since $(b,c)\subset\mathcal{U}$, $f(b)=b$, $f(c)=c$ and $f$ is continuous, there is $d>c$ such that $f((b,c))=(b,d]$.
So, there is $x_1\in (b,c)$ such that $f(x_1)=c$. 
Now, since $f(x_1)>x_1$, there is $x_2\in (b,x_1)$ such that $f(x_2)=x_1$. 
Inductively, we get a strictly decreasing sequence $(x_n)_n\in (b,c)$ such that $f(x_{n+1})=x_n$ for all $n\in \N$. 
It follows that $(x_n)_n$ converges to $b$. 
Let $(n(k))_k\subset\N$ be an increasing sequence such that $\sum_{k=1}^\infty |x_{n(k)}-b|<\infty$. 
Thus, $\mu=\sum_{k=1}^{\infty} \delta(x_{n(k)})-\delta(b)\in \F(M)$.
As we did in Case 3.1, using the Lipschitz map $\varphi(\cdot)= \textup{dist}(\cdot,[a,b])\in \Lip_0(M)$, we get that $\mu \in A_{\widehat{f}}$.
Also, $\bigcup_n\textup{supp}(\widehat{f}^n(\mu))\subset [b,c]$. Thus, Corollary \ref{cor: boundedness of fn(mu)} gives us that $A_{\widehat{f}}$ is dense in $\F(M)$.
The proof is now complete.
\end{proof}

In order to study the last case of Theorem~\ref{theo: Af for R}, i.e. $M=\R$, we need the following elementary lemma.

\begin{lemma}\label{lem: A set and powers}
Let $M$ be a pointed metric space, $f\in \Lip_0(M,M)$ and $k\in\N$. Then $A_{f}=A_{f^k}$ for any $k\geq 1$.
\end{lemma}

\begin{proof}
	Let $k\geq 1$.
	By the very definition of the set $A_f$, it follows that $A_{f}\subset A_{f^k}$. 
	So, we check the reverse inclusion. 
	Let us assume that $A_{f^k}\neq \emptyset$.
	Thus, $\Lip(f)>1$.
	Let $x\in A_{f^k}$ and $R>0$. 
	Reasoning towards a contradiction, assume that $x\notin A_f$. 
	Therefore, there is an increasing sequence $(n(j))_j\subset \N $ such that $d(0,f^{n(j)}(x))< R$. 
	For any $j\in\N$, let $m(j)\in \N$ and $p\in[0,k-1]$ such that  $n(j)=km(j)+p$.
	Observe that $(m(j))_j$ is unbounded.
	Also,
	\[d(0,f^{k(m(j)+1)}(x))= d(f^{k-p}(0),f^{k-p}(f^{n(k)}(x)))\leq \Lip(f)^kR,~\forall j\in\N. \]
	Since $(m(j))_j$ is unbounded, $x\notin A_{f^k}$, which is a contradiction.
\end{proof}

Finally, we can proceed with the proof of Theorem~\ref{theo: Af for R}, that is, the remaining case $M=\R$.

\begin{proof}[Proof of Theorem~\ref{theo: Af for R}]
	
If $M$ is a compact interval, then the conclusion follows from Corollary~\ref{cor: examples}.
If $M=[a,\infty)$ with $a\leq 0$, or $M=(-\infty,a]$ with $a\geq0$, the conclusion follows from Proposition~\ref{prop: semi infinite interval}.
So, let us assume that $M=\R$ with its usual metric.\\

Let $f\in \Lip_0(M,M)$.
Let $\mathcal{U}=\{z>0:~f(z)>z\}\cup\{z<0:~f(z)<z\}$. We split our analysis in three cases.\\

Case 1: $\mathcal{U}\neq \emptyset$. Then, the proof follows the same lines of Case 2 and Case 3 of the proof of Proposition~\ref{prop: semi infinite interval}.\\

Case 2: $\mathcal{U}= \emptyset$ and every orbit generated by $f$ is bounded. The conclusion follows from Proposition~\ref{bounded orbits}.\\

Case 3: $\mathcal{U}= \emptyset$ and there is $x\in \R$ such that $\textup{orb}(f,x)$ is unbounded. Let $n(1)<n(2)<n(3)\in\N$ such that 
\[|f^{n(1)}(x)|<|f^{n(2)}(x)|<|f^{n(3)}(x)|.\]
Thus, there are two different $i,j\in \{1,2,3\}$ such that $\textup{sign}(f^{n(i)})(x)=\textup{sign}(f^{n(j)}(x))$.
Without loss of generality, we assume that $i=1$ and $j=2$ and that $f^{n(1)}(x)>0$. 
Let us fix $y=f^{n(1)}(x)$ and $k=n(2)-n(1)$. 
Then, $f^k(y)>y$ and thus, the set $\mathcal{U}_k= \{z>0:~f^k(z)>z\}\cup\{z<0:~f^k(z)<z\}$ is nonempty. Again, following the lines of the proof of the Case 2 and Case 3 of the proof of Proposition~\ref{prop: semi infinite interval} we get that $A_{\widehat{f^k}}$ is dense in $\F(\R)$. 
Since $\widehat{f^k}=\widehat{f}^k$, Lemma~\ref{lem: A set and powers} gives us that $A_{\widehat{f}}=A_{\widehat{f}^k}$. 
Thus, $A_{\widehat{f}}$ is dense in $\F(\R)$.
The proof is now complete.
\end{proof}



	


\section{The set $\widehat{\Lip}_0(M,N)$}\label{sec: top}
Let $(M,d)$ and $(N,\rho)$ be two complete pointed metric spaces.
In Proposition~\ref{prop: bounded Lip constant} and Remark~\ref{rem: prop bounded lip} we characterized the SOT-convergence of sequences $(\widehat{f}_n)_n$ with respect to properties of the sequence $(f_n)_n\subset \Lip_0(M,N)$. 
In this section we study the set $\widehat{\Lip}_0(M,N)$ with respect to the weak operator topology and strong operator topology inherited from $\mathcal{L}(\F(M),\F(N))$. Recall Definition~\ref{def: sot and wot}.

\begin{proposition}
Let $(f_n)_n\subset\Lip_0(M,N)$ and $g\in \Lip_0(M,N)$ be such that $(\widehat{f}_n)_n$ converges to $\widehat{g}$ with respect to the weak operator topology. Then, $(\widehat{f}_n)_n$ converges to $\widehat{g}$ with respect to the strong operator topology. 
\end{proposition}

\begin{proof}
Let us first show that, for any $x\in M$, we have that $f_n(x)\to g(x)$. 
Indeed, otherwise, there are $x\in M\setminus\{0\}$, a subsequence $(n(j))_k$ and $\varepsilon>0$ such that $\rho (f_{n(j)}(x),g(x))> \varepsilon$.
Let $\varphi\in \Lip_0(N)$ such that $\varphi(g(x))=1$ and $\textup{supp}(\varphi)\subset B(g(x), \varepsilon)$. 
Now, observe that for any $j\in\N$
\[\langle \varphi, \widehat{f}_{n(j)}(\delta(x))\rangle =\langle \varphi, \delta({f}_{n(j)}(x))\rangle =0\neq 1 = \langle \varphi, \widehat{g}(\delta(x))\rangle,\]
which is a contradiction.\\

On the other hand, it is well know that every convergent sequence of bounded linear operators with respect to the weak operator topology is norm bounded. Therefore, there is $C>0$ such that $\|\widehat{f}_n\|<C$ for all $n\in \N$. Thus, $\Lip(f_n)<C$ for all $n\in \N$. 
Now, applying Proposition~\ref{prop: bounded Lip constant} we get that $(\widehat{f}_n)_n$ converges to $\widehat{g}$ with respect to the strong operator topology. 

\end{proof}

\begin{remark}
Note that, in general, the SOT-convergence and WOT-convergence do not coincide sequentially. Indeed, if $F$ denotes the forward shift in $\ell^2$, (with respect to its canonical basis), then $(F^n)_n$ converges to $0$ in the weak operator topology but it does not converge in the strong operator topology. However, these convergences coincide for operators defined on a Banach space with the Schur property. In particular, SOT-convergence and WOT convergence sequentially coincide in $\mathcal{L}(\ell^1)$.
Further information of the Schur property and Lipschitz-free spaces can be found in \cite{PColin}.
\end{remark}

In order to continue our study, we need the following result which is a direct consequence of \cite[Lemma 2.10]{ACP}.

\begin{proposition}\label{weakly closed set}
	Let $(M,d)$ be a complete pointed metric space. Then, $\delta(M)$ is a weakly closed subset of $\F(M)$.
\end{proposition}
\begin{proof}
Let $(x_\alpha)_{\alpha\in\Lambda}\subset M$ be a net such that $(\delta(x_\alpha))_{\alpha\in\Lambda}$ weakly converges to $\mu\in \F(M)$. 
Thanks to \cite[Lemma 2.10]{ACP}, $|\textup{supp}(\mu)|\leq 1$. 
So, $\mu=0$ or there is $\lambda\in\R\setminus\{0\}$ and $x\in M\setminus\{0\}$ such that $\mu=\lambda\delta(x)$.
Assuming that we are in the second case,
by considering functions $\varphi_{\varepsilon}\in\Lip_0(M)$ of the form $\varphi_{\varepsilon}(\cdot)=\textup{dist}(\cdot,M\setminus B(x,\varepsilon))$, we deduce that $\lambda=1$.
\end{proof}

\begin{corollary}
	The set $\widehat{\Lip}_0(M,N)$ is a closed subset of $\mathcal{L}(\F(M),\F(N))$ with respect to the weakly operator topology. Particularly, $\widehat{\Lip}_0(M,N)$ is also a closed set for the strong operator topology.
\end{corollary}
\begin{proof}
	Throughout this proof $\delta_M$ and $\delta_N$ denote the canonical isometry from $M$ to $\F(M)$ and from $N$ to $\F(N)$ respectively.
	Let $(T_\alpha)_{\alpha\in\Lambda}\subset \widehat{\Lip}_0(M,N)$ be a net such that $(T_\alpha)_{\alpha\in\Lambda}$ converges to $T\in \mathcal{L}(\F(M),\F(N))$ in the weak operator topology.
	Let us see that there is $f\in \Lip_0(M,N)$ such that $T=\widehat{f}$.
	Note that, for any $x\in M$ and any $\varphi\in \Lip_0(N)$,
	\[\langle \varphi, T_\alpha(\delta_M(x))\rangle\to \langle \varphi,T(\delta_M(x))\rangle.\]
	Thus, fixing $x\in M$, we have that $(T_\alpha(\delta_M(x)))_{\alpha\in\Lambda}$ converges weakly to $T(\delta_M(x))$. Since $T_\alpha(\delta_M(x))\in \delta_N(N)$ for any $\alpha\in\Lambda$, thanks to Proposition~\ref{weakly closed set}, we get that $T(\delta_M(x))\in \delta_N(N)$.
	Thus, $T(\delta_M(M))\subset \delta_N(N)$. 
	Now, it readily follows that there is $f\in \Lip_0(M,N)$ such that $T=\widehat{f}$.
	Indeed, $f(x)=\delta_N^{-1}(T(\delta_M(x)))$ for any $x\in M$.
\end{proof}
\section{Further comments}

In this paper we have investigated the sets $R_{\widehat{f}}$, $U_{\widehat{f}}$ and $A_{\widehat{f}}$, where $f\in\Lip_0(M,M)$. Our results of Section~\ref{sec: R set} and Section~\ref{sec: A set} suggest that the following two questions hold true, however they are open in their fully generality.

\begin{question}
	For an arbitrary pointed metric space $M$ and $f\in \Lip_0(M,M)$. Is it true that $\textup{int}(R_{\widehat{f}})\neq\emptyset$ implies that $R_{\widehat{f}}=\F(M)$?
\end{question}
\begin{question}
	For an arbitrary pointed metric space $M$ and $f\in \Lip_0(M,M)$. Is it true that $A_{\widehat{f}}\neq\emptyset$ implies that $A_{\widehat{f}}$ is dense in $\F(M)$?
\end{question}
We end this paper with a question motivated by Corollary~\ref{cor: rigid characterizetion}:
\begin{question}
Can we characterize recurrent Lipschitz operators?
\end{question}

\section*{Acknowledgement}
The author is deeply grateful for several fruitful discussion with Frédéric Bayart and Robert Deville. 
The author was partially funded by the Austrian Science Fund (FWF P-36344N).
\bibliographystyle{amsplain}

\end{document}